\documentclass[11pt,a4paper,reqno]{amsart} 

\usepackage{amsmath}
\usepackage{amsthm}
\usepackage{amsfonts}
\usepackage{amssymb}
\usepackage[dvips]{graphicx}
\usepackage{latexsym}

 \DeclareMathOperator{\Ad}{Ad}
 \DeclareMathOperator{\Hom}{Hom}

\newcommand{\ben}{\begin{equation*}}
\newcommand{\een}{\end{equation*}}
\newcommand{\bal}{\begin{aligned}}
\newcommand{\eal}{\end{aligned}}
\newcommand{\baln}{\begin{align}}
\newcommand{\ealn}{\end{align}}
\newcommand{\bma}{\begin{pmatrix}}
\newcommand{\ema}{\end{pmatrix}}

\newcommand{\w}{{\scriptstyle\wedge}\,}
\newcommand{\semi}[1]{\oplus_{#1}}
\newcommand{\lspan}{{\rm span}}
\newcommand{\der}{{d}}
\newcommand{\D}{\mathcal{D}}
\newcommand{\Der}{\mathfrak{D}}
\newcommand{\g}{\mathfrak{g}}
\newcommand{\glg}{GL(2,\real)}
\newcommand{\gl}{\mathfrak{gl}(2,\real)}
\newcommand{\h}{\mathfrak{h}}
\renewcommand{\H}{\mathcal{H}}
\newcommand{\m}{\mathfrak{m}}
\newcommand{\real}{\mathbb{R}}

\newcommand{\V}{{\mathcal V}}
\newcommand{\W}{{\mathcal W}}

\newcommand{\C}{\mathbb{C}}
\newcommand{\CP}{\mathbb{CP}}
\newcommand{\RP}{\mathbb{RP}}
\newcommand{\spp}{\mathbb{S}}
\newcommand{\R}{\mathbb{R}}
\newcommand{\Z}{\mathbb{Z}}
\newcommand{\T}{\mathbb{T}}

\def\p{\partial}
\def\ov{\overline}
\def\l{\lambda}

\newcommand{\hook}{{\setlength{\unitlength}{11pt}   
                   \begin{picture}(.833,.8)
                   \put(.15,.08){\line(1,0){.35}}
                   \put(.5,.08){\line(0,1){.5}}
                   \end{picture}}}
\renewcommand{\d}{{d}}
\def\be{\begin{equation}}
\def\ee{\end{equation}}
\newcommand{\koniec}{\begin{flushright}  $\Box $ \end{flushright}}
\newtheorem{theo}{Theorem}[section] 
\newtheorem{prop}[theo]{Proposition}  

\newtheorem{defi}[theo]{Definition}

\begin{document}
\date{31 July 2010}
\title{{\bf $GL(2, \R)$ structures, $G_2$ geometry and twistor theory}}
\author{Maciej Dunajski}
\address{Department of Applied Mathematics and Theoretical Physics\\ 
University of Cambridge\\ Wilberforce Road, Cambridge CB3 0WA\\ UK.}
\email{m.dunajski@damtp.cam.ac.uk}
\author{Micha{\l} Godli{\'n}ski} 
\address{
SISSA\\
Via Beirut 2/4\\
34100 Trieste\\
ITALY}
\email{godlinsk@sissa.it } 

\begin{abstract}{A $GL(2, \R)$ structure on an $(n+1)$--dimensional manifold is a
smooth point-wise identification of tangent vectors with polynomials in two variables homogeneous of degree $n$.
This, for even $n=2k$, defines a conformal structure of signature $(k, k+1)$
by specifying the null vectors to be the polynomials with vanishing quadratic invariant. We focus on the case $n=6$ and show that 
the resulting conformal structure in seven dimensions is compatible with a conformal $G_2$ structure or its non--compact analogue.
If a $GL(2, \R)$ structure arises on a moduli space of rational curves on a surface with
self--intersection number 6, then certain components of the intrinsic 
torsion of the $G_2$ structure vanish. We give examples of simple 7th order ODEs whose solution curves are rational 
and find the corresponding $G_2$ structures. In particular we show 
that Bryant's weak 
$G_2$ holonomy metric on the homology seven-sphere $SO(5)/SO(3)$ 
is the unique weak $G_2$ metric arising from a rational curve.
}
\end{abstract}
\maketitle
\section{Introduction}
Consider the three--dimensional space $M$ of holomorphic parabolas in $\C^2$. Each parabola is of the form
\[
y=ax^2+2bx+c
\]
and $(a, b, c)$ serve as local holomorphic coordinates on $M$. Two parabolas generically intersect
at two points, and we can define a holomorphic conformal structure
 on $M$ by declaring two points $p$ and $\tilde{p}$ to be null separated iff the corresponding parabolas are tangent. The tangency condition is equivalent
to a polynomial equation
\[
v^3x^2+2v^2x+v^1=0
\]
having a double root. Here $(v^1, v^2, v^3)=(\tilde{c}-c, \tilde{b}-b, \tilde{a}-a)$ is the vector
connecting $p$ and $\tilde{p}$. Calculating the discriminant shows that this vector is null
if $(v^2)-v^1v^3=0$. This quadratic condition defines a flat 
conformal structure on $M=\C^3$.

An immediate question is whether this approach can be generalised to curved conformal structures. One answer goes back to W\"unschmann \cite{wun} who worked in the real
category. The parabolas are integral curves
of a third order ODE $y^{'''}=0$. W\"unschmann has found the necessary and sufficient condition for a general third order ODE so that the conformal structure induced on the solution space by the tangency condition is well defined. The question has also been considered in the context of twistor theory \cite{hitchin} where one is not concerned with differential equations but rather with the algebro--geometric properties of rational curves in a complex two--fold.

How about higher dimensions? It turns out that one can define conformal structures
on certain odd--dimensional moduli spaces of rational curves, but the discriminant
(which is not quadratic for higher degree curves) needs to be replaced by another 
invariant. In this paper we shall consider the seven--dimensional case
and answer the following questions
\begin{itemize}
\item Given a seven--dimensional  family of rational curves, can one define
a conformal complexified $G_2$ structure  on the moduli space $M$ on these curves?
Does this structure admit a real form of Riemannian signature?
\item Can one characterise the curves and the corresponding $G_2$ structures in terms
of a 7th order ODE 
\[
y^{(7)}=F(x, y, y', \dots, y^{(6)})
\]
with $M$ as its solution space?
\end{itemize}
The methods employed in the paper form a mixture of `old' and `new'.
To define the conformal structure on $M$ we use the 19th century classical invariant 
theory (formula (\ref{conf_structure}) Section 3),
but the characterisation of curves builds on the Penrosean holomorphic 
twistor methods. The allowed rational curves must (after complexification) have self--intersection number $6$ in some complex two--fold or a normal bundle 
${\mathcal O}(5)\oplus{\mathcal O}(5)$ in a complex contact three--fold. This allows a point--wise identification of tangent vectors in $M$ with sextic homogeneous polynomials in two variables.
Now the invariant theory can be applied to construct a conformal structure, and the associated $G_2$ three--form $\phi$ (formulae (\ref{7Dconf}) and (\ref{three_form}) Section 4). The ODE approach gives a good handle on the local differential geometry on $M$ and allows expressing the components of intrinsic torsion of the $G_2$
structure (as well as the torsion of the associated Cartan connection)  in terms of the contact invariants of the corresponding 
ODE (Theorems \ref{th.g2} and \ref{th.examples} formulated in Section 5 and proved in Section 10). Here we make an
extensive use of the Tanaka--Morimoto theory of normal Cartan's connection (Sections 8 and 9). These methods allow us show that if the component of the intrinsic $G_2$ torsion taking value in the 27--dimensional irreducible representation
$\Lambda^3({\R^7}^*)$ vanishes, then the resulting $G_2$ geometry admits a Riemannian real form
and (up to diffeomorphisms) it is either flat, or is given by Bryant's weak $G_2$ holonomy
\cite{B87} 
on $SO(5)/SO(3)$, or is given by a seven--parameter family of curves
\[
(y+Q(x))^2+P(x)^3=0,
\]
where the polynomials $(Q(x), P(x))$ are the general cubic and quadratic respectively. These curves have degree six, but
we shall find that they are rational and form a complete analytic family. The corresponding 7th order ODE is
\[
y^{(7)}=\frac{21}{5}\frac{y^{(6)}y^{(5)}}{y^{(4)}}-\frac{84}{25}
\frac{(y^{(5)})^3}{(y^{(4)})^2},\qquad\mbox{where}\quad y^{(k)}=\frac{\p^k y}{\p x^k}
\]
and the associated conformal structure is given by 
(\ref{7Dconf}) and (\ref{cusp_conf}). There exists a choice of the conformal factor such that corresponding  $G_2$ structure is closed, i. e.
\[
d\phi=0, \quad d*\phi=\tau\wedge\phi
\]
for some two--form $\tau$ on $M$. 
\vskip5pt
Most calculations in the second half of the paper  were performed using MAPLE. In particular proving Theorem \ref{th.examples} required solving a system of over 600 quadratic equations for components of curvature and torsion of Cartan's normal connection.
The resulting expressions are usually long and unilluminating and we have not included all of them in the manuscript. Readers who want to verify our calculations can obtain
the MAPLE codes from us.

\subsection*{Acknowledgements}
The idea that a $G_2$ structure may exist on a moduli space of rational curves
with self--intersection number six was suggested to one of us (MD) by Simon Salamon in
1996. We wish to thank Robert Bryant, Boris Doubrov, Stefan Ivanov, Pawel Nurowski and Sasha Veselov for useful discussions.
The work of MG  was partially supported by
the grant of the Polish Ministry of Science and Higher Education N201
039 32/2703.

\section{$GL(2, \R)$ structures}
\noindent
\begin{defi} 
A GL$(2, \R)$ structure on a smooth $(n+1)$ dimensional manifold $M$ is a smooth bundle isomorphism
\be
\label{paracon}
TM\cong \spp \odot \spp \odot \cdots \odot \spp =\mbox{S}^{n}{(\spp)},
\ee
where $\spp\rightarrow M$ is a real rank--two vector bundle, and $\odot$
denotes symmetric tensor product.
\end{defi}
The isomorphism (\ref{paracon}) identifies each tangent space $T_{{\bf t}} M$  
with the space of homogeneous 
$n$th order polynomials 
in two variables. The vectors corresponding to polynomials with repeated root of multiplicity $n$ are called maximally null. A hyper-surface in $M$
is maximally null if its normal vector is maximally null. 

In practice the isomorphism (\ref{paracon}) giving rise to a 
$GL(2, \R)$ structure is specified by
a binary quantic with values in $T^*M$
\be
\label{form}
Q(X_1, X_2)=\sum_{i=0}^n {n\choose i}\theta^{i+1}(X_1)^i (X_2)^{n-i}, \qquad 
{n\choose i}=\frac{n(n-1)\dots(n-i+1)}{i!}.
\ee
Here $(X_1, X_2)$ are coordinates on $\R^2$, and the `coefficients'
in the quantic are given by linearly independent one--forms $\theta^1, \theta^2, \dots, 
\theta^{n+1}$ on $M$.
If $V$ is a vector field on $M$, then the corresponding polynomial is given by
$V\hook Q$, where $\hook$ denotes the contraction of a one--form with a vector field.
If $V=\sum_i v^i \theta_i $ is expressed in a basis $\theta_i$ of $TM$ such that $\theta_i\hook \theta^j=\delta_i^j$, the polynomial is
\be
\label{form1}
\sum_{i=0}^n {n\choose i}v^{i+1}(X_1)^i (X_2)^{n-i},
\ee
with the coefficients $v^i$ being smooth functions on $M$. 
\vskip10pt
Consider a general ODE of order $(n+1)$
\be
\label{ode}
\frac{d^{n+1} y}{dx^{n+1}}=F(x, y, y', \dots, y^{(n)}),
\ee 
where $y'=d y/d x$ etc, whose general solution is
of the form $y=Z(x, {\bf t})$ where 
${\bf t}$ are constants of integration.
Assume that the space of solutions to  (\ref{ode}) is
equipped with a $GL(2, \R)$ structure (\ref{paracon}) such that the  two--parameter family of hyper-surfaces 
given by fixing $(x, y)$ are maximally null.
It has been shown in \cite{DT06} that this imposes
conditions on $F$ which are expressed by vanishing
of $(n-1)$ expressions 
\be
\label{wunshmans111}
W_\alpha[F], \qquad \alpha=1, 2, \dots, n-1
\ee
for the ODE (\ref{ode}). Each expression $W_\alpha$ is a polynomial in 
the derivatives of $F$. The simplest of these is the contact invariant 
\begin{eqnarray*}
W_1[F]&=&{\mathcal D}^2 F_{n}-\frac{6}{n+1}F_{n} {\mathcal D}F_{n}+
\frac{4}{(n+1)^2}(F_{n})^3-\frac{6}{n} {\mathcal D} F_{n-1}\\
&&+\frac{12}{n(n+1)}
F_nF_{n-1}+\frac{12}{n(n-1)}F_{n-2},
\end{eqnarray*}
where
\[
F_k=\frac{\p F}{\p y^{(k)}}, \quad\mbox{and}\quad
{\mathcal D} =\frac{\p}{\p x}+\sum_{k=1}^ny^{(k)}\frac{\p}{\p y^{(k-1)}}
+F\frac{\p}{\p y^{(n)}}.
\]
Moreover if $W_1[F]=W_2[F]=\cdots=W_{m-1}[F]=0$ then $W_m[F]$ is a 
contact invariant of the ODE (\ref{ode}).
The explicit expressions for $W_\alpha$ are unilluminating, but  for 
completeness  we  list the five invariants (in the form given in \cite{Nur})
of 7th order ODEs in Appendix A.

The same invariants have also arisen in other related contexts \cite{Dou,Dou2,Nur}.
The description given by Doubrov is particularly clear. First note that a linearisation
of the ODE (\ref{ode}) around any of its solutions is a linear homogeneous ODE of the form
\be
\label{linear_ode}
(\delta y)^{(n+1)}=p_{n}(x)\delta y^{(n)}+\dots+p_0(x)\delta y,
\ee
where $p_k=\p F/\p y^{(k)}$ is evaluated at the solution.
\begin{theo}\cite{Dou, Dou2}
The expressions {\em(\ref{wunshmans111})\em} vanish if and only if the linear homogeneous
ODE {\em(\ref{linear_ode})} can be brought to a form $\delta y^{(n+1)}=0$ by a coordinate
transformation $(x, y)\rightarrow (\beta(x), \gamma(x)y)$ for some functions $\beta$ and $\gamma$. Vanishing of {\em(\ref{wunshmans111})} is invariant under the contact transformations of the nonlinear ODE {\em(\ref{ode})}.
\end{theo}
The linear homogeneous ODEs of the form (\ref{linear_ode})
have been studied by Wilczynski \cite{Wilczynski}
who gave explicit conditions for their trivialisability in
terms of the functions $p_k$ and their derivatives.

In the simplest nontrivial case $n=2$ the corresponding invariant was already
known to W\"unschmann \cite{wun}. In the case $n=3$ the invariants have been implicitly constructed by Bryant in his study of exotic holonomy \cite{B91} and developed by Nurowski
\cite{Nur_four}.

One source of ODEs for which these contact invariants vanish comes from twistor theory \cite{B91, DT06}. Let $\mathcal{Y}$ be a complex contact three-fold with an embedded rational Legendrian curve with a normal bundle $N={\mathcal O}(n-1)\oplus{\mathcal O}(n-1)$. The moduli space of such curves
is $(n+1)$ dimensional and carries a natural (complexified) $GL(2, \R)$ structure.

The special case is $\mathcal{Y}=P(T\T)$, where ${\T}$ s  complex two--fold ${\T}$ containing
embedded rational curve $L$ with self--intersection number $n$. Such curve 
has a natural lift $\hat{L}$ to $\mathcal{Y}$, given by
$z\in L\rightarrow (z, \dot{z}\in T_zL)$. The lifted curves are
Legendrian with respect to the canonical contact structure on the
projectivised tangent bundle.
The ODE  whose integral curves are given by holomorphic
deformations of $L$ satisfies the $GL(2, \R)$ conditions.

\section{$GL(2, \R)$ conformal structure}
In this section we shall associate a conformal 
structure to a $GL(2, \R)$ structure. From now we assume that $n=2k$ is even.  We shall first recall some classical theory of invariants \cite{Grace_Young}.
Let $V_n\subset\R[X_1, X_2]$ be the $(n+1)$ dimensional space of homogeneous
polynomials of degree $n$. Consider the linear action of $GL(2, \R)$ on $\R^2$ given by
\[
{\tilde X}_1=\alpha X_1+\beta X_2, \qquad {\tilde X}_2=\gamma X_1+\delta X_2, \qquad
\alpha\delta-\gamma\beta\neq 0.
\]
Given  a binary quantic $Q(X_1, X_2)$ (whose coefficients may be numbers, functions, one--forms, ...)
let $\widetilde{Q}({\tilde X}_1, {\tilde X}_2)$ be a binary quantic such that
\[
\widetilde{Q}({\tilde X}_1, {\tilde X}_2)=Q(X_1, X_2).
\]
This induces an embedding  $GL(2, \R)\subset GL(n+1, \R)$, as
the coefficients ${\bf\tilde{\theta}}=(\tilde{\theta}^1, ..., \tilde{\theta}^{n+1})$
are linear homogeneous functions of the coefficients of $Q$.
Recall that an {\em invariant} of a binary quantic 
is a function $I({\bf \theta})$ 
depending on the coefficients ${\bf \theta}=(\theta^1, \theta^2, ..., 
\theta^{n+1})$ such that
\[
I({\bf\theta})=(\det{A})^w I({\bf{\tilde{\theta}}}),\qquad \mbox{where}\qquad   A=
\left (
\begin{array}{cc}
\alpha&\beta\\
\gamma&\delta
\end{array}
\right )\in GL(2, \R).
\]
The number $w$ is called the weight of the invariant. For example if $n=2$
the discriminant $\theta^1\theta^3-(\theta^2)^2$ is an invariant with weight 
$2$.

One of the classical results of the invariant theory is that all invariants
arise from the {\em transvectants} \cite{Grace_Young}. 
\begin{defi}For any homogeneous polynomials
$Q\in V_n, R\in V_m$ the $p$th transvectant is 
\[
<Q, R>_p=\frac{1}{p!}\sum_{i=0}^p{p\choose i}(-1)^i
\frac{\p ^pQ}{\p (X_1)^{p-i}\p (X_2)^i}\frac{\p ^pR}{\p (X_1)^{i}\p (X_2)^{p-i}}\in V_{n+m-2p}.
\]
\end{defi}
In particular specifying $Q=R$ the successive transvectant operations
reduce to elements of $V_0$ which are invariants. The simplest
of these is \[I_0=<Q, Q>_n.\] It vanishes if $n$ 
is odd, and for even $n=2k$ it has weight $n$ and is 
given by
\be
\label{conf_structure}
I_0=\Big\{ 2\sum_{i=0}^{k-1}(-1)^i 
{2k\choose i}
\theta^{i+1}\theta^{2k+1-i}\Big\}
+
{2k\choose k}(-1)^k (\theta^{k+1})^2. 
\ee
In particular if $I_0$ is evaluated for the binary quantic (\ref{form}) defining the 
$GL(2, \R)$ structure where $\theta^i$ are one--forms, 
then $I_0$ should be regarded as a section
of $S^2(T^*M)$. It is well known that a conformal structure $[g]$ is determined 
by specifying null vectors, i.e. sections  $V\in\Gamma(TM)$ such 
that $g(V, V)=0$ for $g\in [g]$.  This gives
\begin{prop}
\label{conformal_from_GL2}
A $GL(2, \R)$ structure on a $(2k+1)$--dimensional manifold $M$ induces 
a conformal structure $[g]$ of signature $(k+1, k)$ or $(k, k+1)$.
A vector field is null w.r.t $[g]$ iff 
the corresponding polynomial has $I_0(V)=0$.
\end{prop}
{\bf Proof.}
The `nullness' of a vector is a quadratic condition and thus
leads to a quadratic bilinear form up to scale. Let a vector $V$ 
correspond to a polynomial  (\ref{form1}). The condition
$I_0(V)=0$,
where $I_0$ is given by (\ref{conf_structure}) is indeed a
quadratic and leads to a symmetric bilinear form
$g(X, Y)=<X, Y>_n$ of signature $(k+1, k)$ or $(k, k+1)$.
\koniec
In general conformal structures induced by $GL(2, \R)$ structures 
form a subclass of all conformal structures of signature $(k+1, k)$, 
except when $k=1$ in which case the two notions are equivalent. For 
$n$  odd (that is for even--dimensional $M$) the bilinear form (compare formula (\ref{conf_structure})) resulting from this definition is anti--symmetric, so does not lead to conformal structures.
\vskip5pt
{\bf Example.} In three dimensions $GL(2, \R)$ structures are the same as conformal structures of Lorentzian signature. This is related  to the isomorphism
\[
SL(2, \R)/\Z_2\cong SO(2, 1)
\]
which underlies the existence of spinors. Let a conformal structure be represented by
a metric $g=\eta_{ij}e^ie^j$ where $\eta=\mbox{diag}(1, -1, -1)$ and $e^i, i=1, 2, 3$ is an orthonormal basis of one--forms. The $GL(2, \R)$ structure is defined by  (\ref{form})
with $\theta^1=e^1+e^3, \theta^2=e^2, \theta^3=e^1-e^3$. A vector $V=v^i\theta_i$ corresponds to a polynomial
\[
v^1+2x v^2+x^2 v^3
\]
where $x=X_2/X_1$. The nullness condition 
\[
I_0(V)=v^1v^3-(v^2)^2=0
\]
is given by vanishing of the discriminant. Thus a vector is null iff the corresponding polynomial has a repeated root. In the standard approach to spinors in three dimensions
one represents a vector by a symmetric two-by-two matrix $V^{AB}$ where $A, B=1, 2$, such that
$g(V, V)=\det{(V^{AB})}$. The non--zero null vectors correspond to matrices with vanishing determinant, which therefore must have rank one. Any such matrix is of the form
$V^{AB}=p^Ap^B$. In our approach the matrix $V^{AB}$ gives rise to a homogeneous polynomial
$V^{AB}X_AX_B$ which, in case of null vectors, has a repeated root $x=-p^1/p^2$.
\vskip5pt
{\bf Example.} The five--dimensional $GL(2, \R)$ structures correspond to special conformal structure in signature $(3, 2)$. The nullness condition can also be described geometrically in this case and the following interpretation is well known in the context of classical invariant theory \cite{Grace_Young, PR}. In the five--dimensional case vectors correspond to binary quartics. A generic quartic will have four distinct roots, and the nullness
condition $I_0(V)=0$ implies that their cross  ratio  is a cube root of unity.
This is the equianharmonic condition. The roots of the quartic, when viewed as points
on the Riemann sphere, can in this case be transformed into vertices of a regular
tetrahedron by M\"obius transformation. Riemannian analogues of such geometries have been studied in \cite{Nur_Bob}.

\vskip5pt
We have been unable to find a geometric interpretation of the null condition
$I_0(V)=0$ in the case of  seven--dimensional $GL(2, \R)$ conformal structures which 
will play a role in the rest of the paper. The vectors correspond to binary
sextics which generically admit six distinct roots $z_1, z_2, \dots, z_6$. In this case
one can also form an $SL(2, \C)$ invariant multi cross-ratio
\[
\frac{(z_1-z_2)(z_3-z_4)(z_5-z_6)}{(z_2-z_3)(z_4-z_5)(z_6-z_1)}.
\]
Let $z_1, z_2, ...,  z_6$ denote positions of six points on a plane.
Given a triangle with vertices $(z_1, z_3, z_5)$, and three points $(z_2, z_4, z_6)$
on the lines $(z_3z_5), (z_5z_1)$ and $(z_1z_3)$ respectively, the lines
$(z_1z_2), (z_3z_4)$ and $(z_5z_6)$ are concurrent iff the multi cross-ratio
is equal to $1$. This is the Ceva theorem.

The theorem of Menelaus states that the points $(z_2, z_4, z_6)$ are colinear if
the multi cross-ratio is equal to $-1$.

We have expressed the invariant $I_0$ in terms
of the roots, hoping to characterise its vanishing it by the Menelaus/Ceva conditions,
but found that the invariant does not vanish in neither of these two cases.
\subsection{Twistor theory}
If the $GL(2, \R)$ structure comes from an ODE, then
induced conformal structure (\ref{conf_structure}) arises from the twistor correspondence described at the end
of Section 2. Here we shall concentrate on the special case when
the Legendrian curves on a complex three--fold are lifts of rational curves
from a two--fold.

Let 
\[
x\longrightarrow (x, y=Z(x, t_1, t_2, \cdots,t_{2k+1})
\]
be a graph of a rational curve $L$ in a complex surface
$\T$ with a normal bundle $N(L)={\mathcal O}(2 k)$. The cohomological obstruction group
$H^1(L, N(L))$ vanishes, and therefore
the Kodaira theorems \cite{kodaira} imply that
the curve belongs to a $2k+1$ dimensional complete
family $\{L_t, {\bf t}\in M\}$ parametrised by  points in a $(2k+1)$--dimensional complex
manifold (the space of solutions to (\ref{ode})). Moreover there exists a canonical isomorphism
\[
T_tM\cong H^0(L_t, N(L_t))
\]
which associates a tangent vector at ${\bf t}\in M$ to a global holomorphic section of a normal bundle $N(L_t)={\mathcal O}(2k)$. Such sections are given by homogeneous polynomials
of degree $2k$ which establishes the existence of a $GL(2, \R)$ structure.

The curve $L_t$ has
self-intersection number $2k$, i.e. 
\[
\delta y=\frac{\p Z}{\p {\bf t}}\delta {\bf t}
\]
vanishes at the zeros of a  polynomial of degree $2k$ in $x=X_2/X_1$. 
In its homogeneous form this polynomial
is a binary form (\ref{form1}) with coefficients $v^\alpha, \alpha=1, \dots, 2k+1$ which 
depend on
$t_\alpha$ and are linear in $\delta t_\alpha$.
A vector at a point in ${\bf t}\in M$ corresponds to a normal vector field to the rational curve $L_t$, i.e. a section
of $N(L_t)={\mathcal O}(2k)$ which is the same as a 
homogeneous polynomial of degree $2k$. The corresponding invariant $I_0$ gives a quadratic
form on $M$ up to a multiple and its vanishing selects the 
null vectors. This determines the conformal structure.
\vskip5pt
In practice one proceeds as follows: If the rational curve is given by
\[
F(x, y, t_\alpha)=0
\]
and its rational parametrisation is
\[
x=p(\lambda, t_\alpha), \qquad y=q(\lambda, t_\alpha)
\]
where $p, q$ are functions rational in  $\lambda\in \CP^1$, then
the polynomial in $\lambda$ giving rise to a null vector is given by the 
polynomial part of
\be
\label{algoritm}
\sum_{\alpha}\frac{\p F}{\p t_\alpha}|_{\{x=p, y=q\}}\delta t_\alpha.
\ee
\section{$G_2$ structures from $GL(2, \R)$ conformal structures}
We shall now restrict to the case $n=6$, and demonstrate that the
seven--dimensional $GL(2, \R)$ manifolds
admit a conformal structure
with a compatible $G_2$ structure. If the associated three--form is closed and co-closed, then
the conformal structure is necessarily flat. We will however find examples of non--trivial 
$G_2$ structures where some components of the torsion vanish. In particular there
is a non--trivial example of weak $G_2$ holonomy compatible with the $GL(2, \R)$ structure.
This example is originally due to Bryant \cite{B87}. In Theorem \ref{th.examples}
we shall show that this example is essentially unique.
\vskip5pt
Consider a $GL(2, \R)$ structure given by the binary form 
\[
Q(x)=\theta^1x^6+6\theta^2x^5+15\theta^3x^4+20\theta^4x^3+15\theta^5x^2+6\theta^6x+\theta^7,
\]
with the corresponding quadratic  invariant (conformal structure) (\ref{conf_structure})
\be
\label{7Dconf}
I_0=\theta^1\theta^7-6\theta^2\theta^6+15\theta^3\theta^5-10(\theta^4)^2.
\ee
Here $x=X_2/X_1$ is an inhomogeneous coordinate on the projective line $\RP^1$.

Use a combination of transvectants to construct a three-form\footnote{Some readers may prefer the two component spinor notation \cite{PR}.
The capital letter indices  $A,B, ...$ take values $1,2$. They are raised and
lowered by a symplectic form represented by an anti-symmetric matrix  
$\varepsilon_{AB}$ on $\R^2$
such that $\varepsilon_{12}=1$. The homogeneous polynomials
are of the form  $Q=Q_{AB...C}\pi^A\pi^B...\pi^C$, where $\pi^A=(X_1, X_2)$.
Then $<Q, P>_{n}=Q_{AB...C}P^{AB...C}$. 
The conformal structure and the three form are given by
\begin{eqnarray*}
I_0&=&  e_{ABCDEF}\odot e^{ABCDEF},\\
\phi&=&{e^{ABC}}_{DEF}\wedge {e^{DEF}}_{GHI}\wedge {e^{GHI}}_{ABC},
\end{eqnarray*}
where $e^{ABCDEF}=e^{(ABCDEF)}$ is $\R^7$ valued one-form
such that
\[
e^{111111}=\theta^1, e^{111112}=\theta^2, e^{111122}=\theta^3, e^{111222}=\theta^4, e^{112222}=\theta^5, e^{122222}=\theta^6, e^{222222}=\theta^7.
\]
}
\be
\label{trv3}
\phi(X, Y, Z)=c<<X, Y>_3, Z>_6,
\ee
where $c$ is some constant, which we chose to be $\sqrt{5/2}$.
\begin{prop}
\label{prop_3_form}
The three--form $\phi$ is compatible with the conformal structure $I_0$:
The vector $V$ is null with respect to $I_0$ iff
\be
\label{bryant}
(V\hook \phi)\wedge (V\hook \phi)\wedge \phi=0
\ee
where $(V \hook\phi)(X, Y):=\phi(V, X, Y)$.
\end{prop}
{\bf Proof.}
Consider the conformal structure induced by the vanishing of
(\ref{conf_structure}). Calculating the components of the three--form
$\phi$ given by (\ref{trv3})  and its dual with respect to $I_0$ gives,
for an appropriate choice of $c$
\begin{eqnarray}
\label{three_form}
\phi&=&{\sqrt\frac{5}{2}}(3\,(\theta^2\wedge \theta^3\wedge \theta^7+\theta^1\wedge \theta^5\wedge \theta^6)+
\theta^4\wedge(\theta^1\wedge \theta^7+6\;\theta^2\wedge \theta^6-15\, \theta^3\wedge \theta^5)),\nonumber\\
\ast\phi&=&\frac{3}{4}(-20\, \theta^1\wedge \theta^4\wedge \theta^5\wedge \theta^6+5\, \theta^1\wedge \theta^3\wedge \theta^5\wedge \theta^7
-20\, \theta^2\wedge \theta^3\wedge \theta^4\wedge \theta^7\nonumber\\
&&
-2\, \theta^1\wedge \theta^2\wedge \theta^6\wedge \theta^7
+30\, \theta^2\wedge \theta^3\wedge \theta^5\wedge \theta^6).
\end{eqnarray}
This is in fact the non-compact form $G^{split}_2$ of the $G_2$ structure, as these
 forms agree with the more usual orthonormal 
frame formulae (see e.g. \cite{B87})
\begin{eqnarray}
\label{standard_defi}
I_0&=&({e^1})^2+({e^2})^2+({e^3})^2-({e^4})^2-({e^5})^2-({e^6})^2-({e^7})^2,\nonumber\\
\phi&=&e^{123}-e^{145}-e^{167}-e^{246}+e^{257}+e^{347}+e^{356},\\
\ast\phi&=&e^{4567}-e^{2367}-e^{2345}-e^{1357}+e^{1247}+e^{1256}+e^{1346},\nonumber
\end{eqnarray}
provided that
\begin{eqnarray*}
e^1&=&\frac{1}{2}(\theta^1+\theta^7),\quad  e^5=\frac{1}{2}(-\theta^1+\theta^7),\quad  
e^2=\frac{\sqrt{6}}{2}(\theta^2-\theta^6)\\
e^6&=&\frac{\sqrt{6}}{2}(\theta^2+\theta^6),\quad  
e^3=\frac{\sqrt{15}}{2}(\theta^3+\theta^5),\\
e^7&=&\frac{\sqrt{15}}{2}(-\theta^3+\theta^5),
\quad  e^4=\sqrt{10}\theta^4.
\end{eqnarray*}
(here $e^{ijk}=e^i\wedge e^j\wedge e^k$ etc). The condition (\ref{bryant}) can now be verified directly. Conversely, given a three form $\phi$, the conformal structure
defined by (\ref{bryant}) is represented by $I_0$ as shown in \cite{B87}. 
\koniec
We have therefore explicitly demonstrated 
that ${\mathfrak{gl}(2, \C)}$ can be embedded in the 
complexification ${\mathfrak{g}_2}^\C\oplus\C$ of ${\mathfrak{g}_2}\oplus\R$,
or equivalently that ${\mathfrak{sl}(2, \C)}$
can be embedded in ${\mathfrak{g}_2}^\C$. This 
follows more abstractly from a theorem of Morozov \cite{morozov} 
which says that
for any nilpotent element $e$ of a complex semi--simple Lie algebra
${\mathfrak{g}}$  there exist $f, h\in {\mathfrak{g}}$ and a homomorphism 
$\rho:{\mathfrak{sl}(2, \C)}\longrightarrow {\mathfrak{g}}$ such that
$
\rho({\bf e})=e, \rho({\bf f})=f,\rho({\bf h})=h,
$ where ${\bf e}, {\bf f},   {\bf h},  $ is the basis
of ${\mathfrak{sl}(2, \C)}$ such that
\[
[{\bf e}, {\bf f}]={\bf h}, \quad 
[{\bf h}, {\bf e}]=2{\bf e},\quad [{\bf h}, {\bf f}]=-2{\bf f}.
\]
\vskip5pt
\subsection{Fernandez--Gray types}
In this paper we follow the standard terminology of $G$--structures
and define a $G_2$ structure on a seven--dimensional manifold 
to be a reduction of the frame
bundle from $GL(7, \R)$ to $G_2$ (or its non--compact analogue
$G^{split}_2$). This 
structure is represented by a three--form $\phi$ in the open orbit of 
$GL(7, \R)$ in $\Lambda^3(M)$. The three--form induces
a metric \cite{B87} on $M$. If $\phi$ is given by (\ref{standard_defi})
then the metric is given by $I_0$. It is Riemannian if
the forms $e^1, \cdots, e^3$ are real and $e^4, \cdots, e^7$ are imaginary
and has signature $(3, 4)$ if all one--forms are 
real\footnote{Another equivalent definition \cite{FG_paper} is to start form
a Riemannian (respectively signature $(3, 4)$) metric $g$ and define a $G_2$ structure
to be a cross product $P:T_tM\times T_tM\rightarrow T_tM$ on each tangent space
which varies smoothly with $t\in M$ and such that $P$ is a bilinear map
satisfying
\[
g(P(X, Y), X)=0, \quad |P(X, Y)|^2=|X|^2|Y|^2-g(X, Y)^2, \quad \forall X, Y,
\]
where $|X|^2=g(X, X)$.
One then {\it defines} the associated three form by
\[
\phi(X, Y, Z)=g(P(X, Y), Z).
\]
This cross product equips each tangent space with the algebraic structure
of pure octonions (respectively pure split octonions).}. 
The latter case corresponds to the non--compact form $G^{split}_2$. Proposition 
\ref{prop_3_form}
shows that seven--dimensional $GL(2, \R)$ structure is equivalent
to a further reduction of the frame bundle from $\R^+\times G_2\subset \R^{+}\times SO(3, 4)$
to $GL(2, \R)$.

We do not assume anything  about the closure of the three--form $\phi$ or its 
dual. Various types of $G_2$ structures (or their non--compact analogues)
are characterised by  a representation theoretic decomposition of 
$\nabla\phi$, where $\nabla$ is the Levi--Civita connection of the metric 
induced by $\phi$. Following \cite{FG_paper, B03} we have
\begin{eqnarray}
\label{FG_deco}
d\phi&=&\lambda*\phi+\frac{3}{4}\Theta\wedge\phi+*\tau_3\\
d*\phi&=&\Theta\wedge *\phi-\tau_2\wedge\phi,\nonumber
\end{eqnarray}
where $\lambda$ is a scalar, $\Theta$ is a one--form, $\tau_2$ is a two--form
such that $\tau_2\wedge\phi=-*\tau_2$ and $\tau_3$ is a three--form 
such that $\tau_3\wedge\phi=\tau_3\wedge *\phi=0$.
The forms $(\lambda, \Theta, \tau_2, \tau_3)$ can be 
interpreted as components of 
intrinsic torsion of a natural connection of $G_2$ structure.
To define this connection apply the canonical decomposition
$\mathfrak{so}(7)=\mathfrak{g}_2\oplus \R^7$ (or its $\mathfrak{so}(3, 4)$
analog) to the Levi--Civita connection. If $\gamma$ is the 
$\mathfrak{so}(7)$--valued connection one--form, then writing
$\gamma=\hat{\gamma}+\tau$ defines a connection  with
torsion (not to be confused with the torsion $T$ of the $\gl$--valued connection $\Gamma$ studied in Sections 7 and 10) represented by a one--form $\hat{\gamma}$ with values in $\mathfrak{g}_2$.
See e.g. \cite{B03} for details. If the three--form is only defined up to a multiple
by a non--zero function (as it is the case in this paper) then 
$\lambda, \tau_2 $ and $\tau_3$ scale with appropriate weights and $\Theta$ 
transforms like a Maxwell field. More precisely conformal rescalling $g\rightarrow e^{2f}g$
leaves (\ref{FG_deco}) invariant if
\be
\label{FG_scall}
\phi\rightarrow e^{3f}\phi, \quad \lambda\rightarrow e^{-f}\lambda, \quad
\Theta\rightarrow \Theta+4df, \quad \tau_2\rightarrow e^{f}\tau_2, 
\quad \tau_3\rightarrow e^{2f}\tau_3.
\ee
If all components of the torsion vanish, then the $G_2$ structure gives
rise to $G_2$ holonomy and the resulting metric is Ricci--flat. Such $G_2$ structures are sometimes called {\it integrable} or more correctly
{\it torsion--free}. If $\Theta=\tau_2=\tau_3=0$ then
the metric is Einstein with non--zero Ricci scalar and one
speaks of {\it weak $G_2$ holonomy}. 
If $\lambda=\Theta=\tau_3=0$ then the $G_2$ 
structure is {\it closed} (see e.g. \cite{ivanov}). 

The representation theoretic decomposition of the torsion is as follows:
\begin{itemize}
\item $\lambda$ is a function and $\lambda\,\phi$ belongs to the
1-dimensional irreducible representation
$\W_1\subset\Lambda^3\real^{*7}$ of $G_2$.

\item The 2-form $\tau_2$ belongs to the 14-dimensional irreducible
representation $\W_2\subset \Lambda^2\real^{*7}$.
\item The 3-form $\tau_3$ belongs to the 27-dimensional irreducible
representation $\W_3\subset\Lambda^3\real^{*7}$.
\item The \emph{Lee 1-form} $\Theta$ belongs to the 7-dimensional
representation $\W_4=\real^{7*}$.
\end{itemize}
Equations \eqref{FG_deco} uniquely define $\lambda$, $\tau_2$,
$\tau_3$ and $\Theta$. Vanishing of these objects defines the
Fernandez--Gray $\W$ type of $G_2$ geometry: if none of them
vanishes the geometry is of generic type $\W_1+\W_2+\W_3+\W_4$, if
$\lambda=0$ then the geometry is of type $\W_2+\W_3+\W_4$, when
$\tau_2=0$ we have the type $\W_1+\W_3+\W_4$ and so on. There are
sixteen $\W$ types.

Proposition \ref{prop_3_form} demonstrates that
the seven-dimensional $\glg$ geometry 
is a special case of conformal split $G_2$ geometry
as $\glg \subset \real\times  G^{split}_2$. 
The representations of $G^{split}_2$ decompose into
irreducible representations of $\glg$ as follows.
\be\label{e.wroz}
\begin{aligned}
&\W_1=V^1, \\
&\W_2=V^3\oplus V^{11},  \\
&\W_3=V^5\oplus V^9 \oplus V^{13},\\
&\W_4=V^7,
\end{aligned}
\ee 
where $V^k$ is the $k$--dimensional representation space 
$\mbox{S}^{k-1}{(\spp)}$.
Hence $\tau_2$ and $\tau_3$ have a priori two and three
irreducible components under action of $\glg$.

A $\glg$ geometry defines a whole conformal class of
$ G^{split}_2$ geometries, hence we may only talk about those $\W$
types which are invariant with respect to conformal rescalings 
(\ref{FG_scall}). In particular vanishing
of $\Theta$ is not conformaly invariant.
However, $\der\Theta$ is a well defined 2-form, in particular the
condition $\der\Theta=0$ means that in some conformal gauge
$\Theta$ vanishes locally. In general $\der\Theta$ is a 2-form
decomposing according to \be\Lambda^2\real^7=V^3\oplus V^7\oplus
V^{11}. \ee
\section{$G_2$ structures from ODEs}
Now we are ready to give the relations between the
the intrinsic torsion of the 
split $G_2$ structure and the contact invariants of the 7th order ODE
\be 
\label{7th_ODE}
y^{(7)}=F(x, y, y', \dots, y^{(6)}). 
\ee
In the following theorems (which will be established in Section \ref{sec.gl2})
we shall assume the vanishing of the conditions $W_\alpha$ (Appendix A)
which are necessary and sufficient for an ODE to give rise to a $GL(2, \R)$ geometry.
\begin{theo}\label{th.g2}
Let the 7th order ODE {\em(\ref{7th_ODE})}
admit the $\glg$ geometry on the solution space. The following conditions hold
\begin{align*}
&\begin{gathered}\text{no $\W_1$ component}
\\ (\lambda= 0)
\end{gathered} &&\iff &&
\begin{aligned} &F_{66}(9\D
F_{6}-\tfrac97F_6^2-15F_5)+\\&12F_{65}F_6+14F_{55}-\tfrac{84}{5}F_{64}=0.\end{aligned}
\\&&&&&&&&&\\
&\begin{gathered} \text{no $\W_2$ component}
\\ (\tau_2= 0) \end{gathered} &&\iff
 &&21 \D F_{66} + 14 F_{65} + 15 F_6 F_{66}=0. \\&&&&&&&&&\\
&\begin{gathered} \text{no $\W_3$ component}
\\ (\tau_3= 0) \end{gathered} &&\iff &&F_{66}=0.
\end{align*}

The 2-form $\der\Theta$ falls into components in irreducible representations $V^3$, $V^7$ and
$V^{11}$. The $V^3$-part is expressed algebraically by $\lambda$,
$\tau_2$ and $\tau_3$. In particular it vanishes if $\tau_2$
vanishes. 
The $V^7$-part of $\der\Theta$ vanishes iff \ben \bal
&(\D F)_{66}F_{66}+\frac32(\D
F)_6F_{666}-\frac{12}{7}F_{666}F_6^2-4F_{666}F_5\\
&+2F_{665}F_6-\frac{14}{5}F_{664}+\frac73F_{655}-
\frac43F_{66}F_{65}-\frac{16}{7}F_{66}^2F_6=0.\eal \een 
The
$V^{11}$-part of $\der\Theta$ vanishes iff \ben F_{666}=0. \een
\end{theo}
The non-generic $\W$ types are characterised by the following result 
\begin{theo}\label{th.examples}
There are only three conformal split $G_2$ geometries from ODEs of
type $\W_1+\W_2+\W_4$.
\begin{itemize}
\item[1.] The flat geometry of $y^{(7)}=0$, which is the only case
admitting holonomy $ G^{split}_2$. \item[2.] The geometry of \be
\label{submax_ODE}
y^{(7)}=7\frac{y^{(6)}y^{(4)}}{y^{(3)}}+\frac{49}{10}\frac{(y^{(5)})^2}{y^{(3)}}
-28\frac{y^{(5)}(y^{(4)})^2}{(y^{(3)})^2}+\frac{35}{2}
\frac{(y^{(4)})^4}{(y^{(3)})^3}. \ee
This is the only geometry of type $\W_1+\W_4$. The Lee form is
closed, so that in certain conformal gauge it is the
nearly-parallel ($\W_1$) geometry of $SO(3,2)/SO(2,1)$. \item[3.]
The geometry of \be
\label{cusp_ODE}
y^{(7)}=\frac{21}{5}\frac{y^{(6)}y^{(5)}}{y^{(4)}}-\frac{84}{25}
\frac{(y^{(5)})^3}{(y^{(4)})^2},
\ee which is of type $\W_2+\W_4$. The Lee form is closed, so in
certain conformal gauge it is a closed $G_2$ structure ($\W_2$).
\end{itemize}
\end{theo}
The $G_2$ geometry associated with (\ref{submax_ODE}) has two real forms:
The homogeneous space 
$SO(3, 2)/SO(2, 1)$ which yields a weak $G_2$ metric in signature $(3, 4)$ and
$SO(5)/SO(3)$ which gives a Riemannian metric. The later metric was first constructed by
Bryant in his seminal paper \cite{B87} without using the ODE or twistor techniques.
Theorem \ref{th.examples} implies that up to diffeomorphisms of $M$ this is the only weak
$G_2$ metric arising from an ODE. The ODE (\ref{submax_ODE}) has appeared is several other contexts\footnote{{\em Note added in February 2012.} 
The general solution to this ODE has been constructed in \cite{DS10}.}
See \cite{sokolov, olver}.
\vskip5pt
The ODE (\ref{cusp_ODE}) has an elementary solution given by certain rational curve
which will be analysed in the next section. The solution to
(\ref{submax_ODE}) can also be constructed in terms of rational curves, but explicit description in this case is more involved \cite{DS10}.
\vskip5pt
We note that there exists at least one more connection between differential equations and non--compact $G_2$: The holonomy of an ambient metric
associated to Nurowski's $(3, 2)$ conformal structure \cite{Nurowski} 
is contained in $G^{split}_2$.

\section{Examples}
In this section we give some examples.
The first three
arise on moduli spaces of rational curves by 
twistor theoretic techniques. The last one comes from an
ODE satisfying the W\"unshmann conditions (\ref{wunshmans111}).

We shall write the general 7th order ODE 
(\ref{7th_ODE})
as
\[
y^{(7)}=F(x, y, p, q, r, s, t, u),
\]
where $p=y', q=y'', r=y^{(3)}, s=y^{(4)}, t=y^{(5)}, u=y^{(6)}$.

The examples below can be partially classified by the dimension of the group 
of contact symmetries (recall that the maximal symmetry group
of the trivial ODE $y^{(7)}=0$ is eleven dimensional and given by $GL(2, \R)\ltimes \R^7$).

\subsubsection{Example 1.} Consider a hyperelliptic curve of degree 6 with 2 cusps. The general sextic  has genus 10 and so is not rational, but in our case the genus is zero and a rational parametrisation exist. To see it write the curve as
\be
\label{cuspic_sextic}
(y+Q(x))^2+P(x)^3=0,
\ee
where $(Q, P)$ are general cubic and quadratic respectively which we write as
\[
Q(x)=q_0+q_1x+q_2x^2+q_3x^3, \quad P(x)=p_3(x-p_2)(x-p_1).
\]
This has three singular points. Two double points at $(p_1, -Q(p_1))$ and 
$(p_2, -Q(p_2))$ of type $[2, 1, 1]$ (see e. g. \cite{alg_geom_book}) and one point
of order 4 at infinity, of type $[4, 8, 2]$ which can be seen by writing 
(\ref{cuspic_sextic}) in the homogeneous coordinates. Calculating the genus yields
\[
{\tt g}=\frac{5\cdot4}{2}-1-1-8=0,
\]
as the quadruple point at infinity is not ordinary and has the $\delta$--invariant equal to 
$8$. The rational parametrisation can now be found 
\begin{eqnarray*}
x(\lambda)&=&\frac{p_1+p_2 \lambda^2}{\lambda^2+1},\\
y(\lambda)&=&{p_3}^{3/2}(p_1-p_2)^3\frac{\lambda^3}{(\lambda^2+1)^3}-Q(x(\lambda)).
\end{eqnarray*}
Eliminating the parameters $(p_1, p_2, p_3, q_0, \dots, q_3)$ between
(\ref{cuspic_sextic}) and its six derivatives yields the 7th order 
ODE characterising
the sextic (\ref{cuspic_sextic})
\[
\frac{d^7y}{d x^7}=\frac{21}{5}\frac{ut}{s}-\frac{84}{25}\frac{t^3}{s^2},
\]
which is the ODE (\ref{cusp_ODE}) from Theorem \ref{th.examples}.

Using the prescription (\ref{algoritm}) we find that the 
conformal structure and the associated three--form are represented by
(\ref{7Dconf}) and (\ref{three_form}) with
\begin{eqnarray}
\label{cusp_conf}
\theta^1&=&-2\Omega\sum_{\alpha=0}^3 (p_2)^\alpha dq_\alpha, \quad 
\theta^7=-2\Omega\sum_{\alpha=0}^3 (p_1)^\alpha dq_\alpha,\\
\theta^2&=&-\frac{\Omega}{2}(p_2-p_1)^2(p_3)^{3/2} dp_2, \quad 
\theta^6=\frac{\Omega}{2}(p_2-p_1)^2(p_3)^{3/2} dp_1,\nonumber\\
\theta^3&=&-\frac{\Omega}{15}(3dq_0+(2p_2+p_1)dq_1+(2p_1p_2+(p_2)^2 )dq_2
+3p_1(p_2)^2dq_3),\nonumber\\
\theta^5&=&-\frac{\Omega}{15}(3dq_0+(2p_1+p_2)dq_1+(2p_1p_2+(p_1)^2 )dq_2
+3p_2(p_1)^2dq_3),\nonumber\\
\theta^4&=&-\frac{3\Omega}{20}(p_2-p_1)^2 \sqrt{p_3}\;d(p_3(p_2-p_1)),\nonumber
\end{eqnarray}
where $\Omega=(p_1-p_2)^{-12/5}{(p_3)}^{-9/10}$.

This conformal $G_2$ structure can be analytically continued to Riemannian 
signature:
Setting $p_2=p, p_1=\ov{p}$ where $p\in \C$ and keeping
$(q_0, q_1, q_2, q_3, p_3)$ real gives purely imaginary $\theta^4$ and
\[
\theta^7=\ov{\theta^1}, \qquad \theta^6=-\ov{\theta^2}, \qquad \theta^3=\ov{\theta^5}.
\]
The corresponding conformal structure 
is positive definite and the three form $\phi$ is real. It gives rise to
a closed (in a sense of decomposition (\ref{FG_deco})) $G_2$ structure
as $d\phi=0, d*\phi=-\tau_2\wedge\phi$ in agreement with Theorem \ref{th.examples}.
This Theorem also implies that up to diffeomorphisms this is the only closed $G_2$
structure arising from ODEs\footnote{{\em Note added in February 2012.}
There is also a co-closed example \cite{DD} which arises on the parameter space
of cuspidal cubics in $\CP^2$. The corresponding 7th order ODE goes back
at least to Wilczynski \cite{Wilczynski}.}.

\subsubsection{Example 2.} Consider a rational curve in $\CP^1\times \CP^1$ of bidegree $(1, k)$
\[
y=\frac{r_0+r_1x+ \dots +r_k x^k}
{s_0+s_1x+\dots+ s_k x^k}.
\]
It has self-intersection number $2k$, and the
enumerator of the perturbed curve (section of a normal bundle) 
$\delta y$ defines the conformal structure (\ref{conf_structure}) 
with
\[
\theta^{i+1}= {2k\choose i}^{-1} \sum_{\alpha+\beta=i}
(r_\alpha d s_\beta-s_\beta d r_\alpha),\qquad i=0, \dots, 2k.
\]
This conformal structure is defined on a hypersurface where the resultant
of the denominator and enumerator in $y$ has a non--zero fixed value.
Alternatively we can fix the ambiguity by choosing affine coordinates, 
say $r_k=1$. Now restrict to the seven--dimensional case  $k=3$.
This also gives 
$\phi\wedge d\phi=0$.
The corresponding $7$th order ODE is
\[
\frac{d^7 y}{d x^7}=\frac{P}{Q},
\]
where
\begin{eqnarray*}
P&=&420q^{2}u^{2}+2520qst^{2}-1680qrut-2100qs^{2}u
-504pt^{3}\\
&&+1680r^{2}t^{2}-6300\,trs^{2}+840tups+2625s^{4}-280u^{2}rp+2800ur^2s\\
Q&=&360q^2t-1200rqs-240rtp
+800r^3+300s^2p.
\end{eqnarray*}
This example has six--dimensional group of point symmetries,
given by the M\"obius transformations of $x$ and $y$.
\subsubsection{Example 3.}
We can construct less trivial conformal structures and the 
associated three forms by generalising the last example, and
taking a double covering  of a neighbourhood of a
non--singular curve of bidegree $(1, 6)$ branched along a fixed curve.
Consider a $(1, 6)$ curve in $\CP^1\times \CP^1$
\be
\label{bidegree1}
y=\frac{R(x)}{S(x)}
\ee
where
\[
S=s_0+s_1x+...+s_6x^6, \qquad R=r_0+r_1x+...+r_6x^6.
\]
This curve has normal bundle ${\mathcal O}(12)$, and is parametrised by $\CP^{13}$ minus 
a hypersurface where both polynomials have common factor. We take the branch locus
to be  the $(1, 6)$ curve
\[
y=x^6.
\]
The curves in a covering space we are constructing project to those curves
(\ref{bidegree1}) which meet the branch locus in seven points to second order.
Thus 
\[
x^6S(x)-R(x)=(t_0+t_1x+...+t_6x^6)^2.
\]
This gives 13 conditions on 20 coefficients $(s, r, t)$, leaving the seven dimensional moduli space of curves.
\subsubsection{Example 4.} In \cite{DT06} it was shown that
the moduli space of solutions to the 
ODE \[
\frac{d^{n+1} y}{d x^{n+1}}= \Big(\frac{d^{n} y}{d x^{n}}\Big)^{\frac{n+1}{n}}
\]
admits the $GL(2, \R)$ structure.
Consider a solution curve $x\rightarrow (x, y(x))$ and its perturbation $\delta y$ 
\begin{eqnarray*}
y&=&t_1+t_2 x+...+t_{n}x^{n-1}-\frac{n^n}{(n-1)!}\ln(x+t_{n+1})\\
\delta y&=&\frac{1}{x+t_{n+1}}\Big( \Big(-\frac{n^n}{(n-1)!}\delta t_{n+1}+t_{n+1}\delta t_1\Big)+ \sum_{i=1}^{n-1}(\delta t_i+t_{n+1}\,\delta t_{i+1})\,x^i+
\delta t_{n}\, x^{n}\Big).
\end{eqnarray*}
The enumerator of the polynomial $\delta y$ defines a conformal structure
(\ref{conf_structure}) with
\[
\theta^1=-\frac{n^n}{(n-1)!}dt_{n+1}+t_{n+1}dt_1, \quad
\theta^{i+1}={2k\choose i}^{-1}\Big(d t_i+t_{n+1}\,d t_{i+1}\Big),  
\quad \theta^{2k+1}=d t_{2k},
\]
where $i=1, \cdots, 2k-1$.
We can now specify $2k=6$ and construct the three--form. We find that
$\phi\wedge\d \phi=0$ so that $\lambda=0$ but there is no conformal scale which makes
$\phi$ closed.
\section{Construction of the Cartan Connection}
We shall now describe the $\glg$ and conformal $G^{split}_2$ structures arising from
 7th order ODEs by constructing a $\gl$-valued linear connection on $M$. The basic object in this description is the torsion, which contains lowest order invariants of the $\glg$ geometry, identifies the Fernandez--Gray types of the associated conformal ${G^{split}_2}$ geometry, and expresses these quantities by contact invariants of the underlying ODE.
This approach will give us a better handle  on the various torsion components.
Our aim is to express these components in terms of invariants of the 7th
order ODE and eventually prove Theorems \ref{th.g2} and \ref{th.examples}
Our treatment of Cartan's connection follows closely that of \cite{Nur}.

We shall use an equivalent form of  Definition \ref{paracon}
and regard a $\glg$ geometry on a manifold $M$ as a reduction of the frame bundle $F M$ to
a $\glg$-sub-bundle, where $\glg\subset GL(n,\real)$ acts
irreducibly in each tangent space \cite{Nur}.
We shall focus on the case $n=6$ where 
\[
\glg\subset\real^+\times G^{split}_2\subset \R^+\times SO(3, 4)
\]
holds (see Proposition \ref{prop_3_form}).
The central role will be played by the 
six-jet space $J^6$ and its description via the Tanaka--Morimoto theory, \cite{Tan2, Mor2}, which is a special version of Cartan's method of equivalence. 
We shall first construct a $\gl\semi{.}\real^7$-valued Cartan connection $\Omega$ on a bundle over $J^6$ and then  re-interpret $\Omega$ from the point of view of the $\glg$ structure. The conditions for the existence of the geometry appear to be certain linear conditions for the curvature of $\Omega$. If they are satisfied, then the $\gl$-part of $\Omega$ is the desired linear connection on $M$.
\subsection{Jet space}
Let us consider the space $J^6$ of six-jets of functions from $\real$ to $\real$.  It is an eight-dimensional real manifold, locally parametrised by $(x,y,y_1,\ldots,y_6$), and such that each curve $x\mapsto (x,f(x))$ in the $xy$-space has a unique lift to a curve in $J^6$ given by $x\mapsto (x,f(x),f'(x),\ldots,f^{(6)}(x))$ in the above coordinate system. This gives a distinguished family of all curves lifted from the $xy$-space in $J^6$. One may encode this family in a coordinate-free language of distributions. Let us fix a point $w\in J^6$ and consider all lifted curves through $w$. The linear span of their tangent vectors at $w$ is a two-dimensional subspace $C_w$ in $T_wJ^6$. The collection $C=\cup_w C_w$ is by definition \emph{the contact distribution} on $J^6$. It is generated by two vector fields
\ben
\Der=\partial_x+y_1\partial_y+y_2\partial_{y_1}+\ldots+y_6\partial_{y_5}\qquad \text{and} \qquad \partial_{y_6}.
\een
Given the distribution $C$ we define
\ben \partial C=[C,C],\quad \partial^2 C=[\partial C,C],\quad \ldots \quad \partial^5 C=[\partial^4 C,C],\quad \partial^6C=TJ^6. \een
The distributions constitute a filtration, that is 
\be\label{e.filt} C\subset \partial C\subset\ldots \subset \partial^5C\subset\partial^6 C=T J^6 \ee
and
\be\label{e.fc} [\partial^i C, \partial^j C]= \partial^{i+j+1} C. \ee

The diffeomorphisms of $J^6$ which preserve $C$ are called contact transformations. The well--known Lie-B\"acklund theorem states 
that all contact transformations of $J^6$ are uniquely defined by the contact transformations of $J^1$ and have the following form 
\be
\label{e.cont}
   x\mapsto\bar{x}(x,y,y_1),\quad
   y\mapsto\bar{y}(x,y,y_1),\quad
   y_1\mapsto\bar{y}_1(x,y,y_1),
\ee   
and for higher order jet coordinates
\[ 
{y}_{k+1}\mapsto \frac{\Der\bar{y}_{k}}{\Der\bar{x}}, \qquad
k=1,2,\ldots, n\quad .
\]
The functions $\bar{x},\bar{y}$ and $\bar{y}_1$ in \eqref{e.cont} are not arbitrary but subject to the condition \ben
\bar{y}_1=\frac{\Der\bar{y}}{\Der\bar{x}}. \een
The contact transformation preserve the whole filtration \eqref{e.filt}.

Now consider the 7th order ODE (\ref{7th_ODE}). Any solution $y=f(x)$ of the equation is uniquely defined by
a choice of $f(x_0)$, $f'(x_0)$, $\ldots$, $f^{(6)}(x_0)$ at some
$x_0$. Since this choice of initial data is equivalent to a choice of a point in
$J^6$ there exists exactly one lifted curve $x \mapsto
(x,f(x),f'(x), \ldots, f^{(6)}(x))$ through any point of $J^6$.
Therefore the solutions form a one-dimensional foliation in $J^6$. The corresponding 
tangent distribution is spanned by 
\ben
\D=\partial_x+y_1\partial_y+y_2\partial_{y_1}+\ldots+y_6\partial_{y_5}+F\partial_{y_6}.
\een
An important consequence of this  is that $J^6\to M$ is locally a line bundle,
where $M$ is the solution space of the 7th order ODE
(\ref{7th_ODE}). 
\begin{defi}
The contact geometry of 7th order ODEs is the jet space $J^6$ equipped with
\begin{itemize}
 \item The filtration $C\subset\ldots\subset\partial^5 C\subset \partial^6 C=T J^6$. 
 \item The foliation by the solutions, tangent to the field $\D$.
\end{itemize}
\end{defi}

One may associate to the contact geometry of the ODEs a sub-bundle 
\be\label{r.p}
\tilde{G}\to\tilde{P}\to J^6 \ee of the frame bundle $FJ^6$: The structure group $\tilde{G}$ is the lower triangular group preserving the filtration and the 1-distribution $\lspan\{\D\}\subset C$ tangent to solutions. 

\subsection{Cartan connection}
The main object we use in the construction is Cartan connection defined here as in \cite{Kob}.
\begin{defi}\label{def.cc}
Let $M$ be a manifold of dimension $n$, $G$ a Lie group, $H$ a closed subgroup of $G$ with $\dim G/H=n$ and $H\to P\xrightarrow{\pi} M$ a principal bundle. A Cartan connection of type $(G,H)$ on $ P$ is a one--form ${\Omega}$
with values in the Lie algebra $\g$ of $G$ satisfying the following conditions:
\begin{itemize}
\item[i)] ${\Omega}_u:T_u P\to\g$ for every $u\in P$ is an
isomorphism of vector spaces. \item[ii)] $\Omega(A^*)=A$
for every $A\in\h$ and the corresponding fundamental field $A^*$.
\item[iii)] $R^*_h{\Omega}=\Ad(h^{-1}){\Omega}$ for every $h\in
H$.
\end{itemize}
\end{defi}

The curvature of a Cartan connection is a $\g$-valued 2-form on
$P$ defined by \ben
{K}(X,Y)=\der{\Omega}(X,Y)+\frac{1}{2}[{\Omega}(X),{\Omega}(Y)].
\een
If $\Omega$ is given in a matrix representation then
\be\label{e.curv}
K=\der\Omega+\Omega\w\Omega.
\ee
The curvature is \emph{horizontal}, that is it vanishes on
each vertical vector field: \be\label{e.hor} {K}(X,\,\cdot\,)=0\quad\text{if}\quad\pi_*(X)=0.\ee 
Horizontality of the curvature is locally
equivalent to the property iii) in Definition \ref{def.cc} . Cartan connections with
vanishing curvature are called flat.

\vskip5pt
We are now in position to describe the construction of $\glg$ geometry on the solution space. We start from the bundle $\tilde{P}$ of the contact geometry of ODEs. We are interested in invariants of this geometry. The filtration  is preserved by the contact transformations but the foliation of solutions is not, and generates the contact invariants of the underlying ODE. However, the situation further is complicated by the fact that the object generating the invariants -- a Cartan connection -- exists on a sub-bundle $P\subset\tilde{P}$ rather than $\tilde{P}$ itself. Using the Tanaka--Morimoto theory we shall construct the sub-bundle $H\to P \to J^6$ together with a Cartan connection $\Omega$ of type $(\glg\rtimes\real^7,H)$, where $H$ isomorphic to the group of triangular $2\times2$ matrices. The curvature $K$ of $\Omega$ contains all the local information about the contact geometry of the ODEs. The contact invariants are either components of $K$ or certain combinations of their derivatives of sufficiently high order.

The jet space $J^6$ is a bundle over the solution space $M$ and $P\to M$ is also a principal bundle with the structure group $\glg$. That $\Omega$ generates the $\glg$ geometry on $M$ only if certain  conditions (which we will determine) hold.
First of all, we ask whether $\Omega$ (which is
a Cartan connection on $P\rightarrow J^6$) satisfies the conditions for the Cartan connection of $P\rightarrow M$. It holds if and only if 
\ben {K}(X,\,\cdot\,)=0\text{ for all }X\text{ vertical with respect to }P\to M. \een 
This condition is not satisfied automatically but only holds for the ODEs
with vanishing W\"unschmann invariants (Appendix A).

The Cartan connection $\Omega$ on $P\to M$ is of type $(\glg\rtimes\real^7,\glg)$. It naturally decomposes into the $\real^7$-part and the $\gl$-part. The former behaves like a canonical form $\theta$ on a principal bundle and turns $P$ into a sub-bundle of the frame bundle $FM$. The latter is a linear $\gl$-valued connection $\Gamma$ on $P$. Together $\theta$ and $\Gamma$ define a $\glg$ geometry on $M$. The torsion $T$ and curvature of $\Gamma$ contain the information about local invariants of the geometry, which are in turn expressed by contact invariants of the underlying ODE, since $\Omega$ also describes the contact geometry of the ODEs.

\vskip5pt

{\bf Example.} For the trivial equation $y^{(7)}=0$,  all the objects may be immediately constructed by means of the symmetry group. The full group of contact symmetries is $\glg\ltimes\real^7$. Its action on $J^6$ is transitive and turns it into a homogeneous space $\glg\ltimes\real^7/H$, where $H$ is isomorphic to the group of triangular $2\times 2$ matrices. Thus we have the bundle $H\to P\to J^6$ and $P=\glg\ltimes\real^7$ locally. The connection ${\Omega}$, flat in this case, is given by the Maurer-Cartan form on $P$.
\section{The Tanaka--Morimoto theory}
We turn to detailed description of the  construction. First of all, we briefly describe the general pattern, next we apply it to our case. 
The references for this subsection are \cite{Tan1, Tan2, Mor1, Mor2} and \cite{DKM}.
The contact geometry of ODEs contains the filtration \eqref{e.filt} which is encoded by the graded tangent bundle $gr\, TJ^6$, denoted here by $gr$ for short. Its fibre over $w\in J^6$ is $gr(w)=\bigoplus_{i=1}^7 gr_{-i}(w)$, where  
 \begin{align*} &gr_{-1}(w)=C_w, & &gr_{-2}(w)=\partial C_w/C_w,\,\ldots \\
 \ldots\,\,&gr_{-6}(w)=\partial^5 C_w/\partial^4 C_w & &gr_{-7}(w)=T_wJ^6/\partial^5 C_w.\end{align*}
The relation \eqref{e.fc} implies that $gr(w)$ carries the structure of a nilpotent graded Lie algebra, that is 
\ben [gr_{-i}(w),gr_{-j}(w)]\subset gr_{-i-j}(w),\quad\text{and } gr(w)\text{ is generated by } gr_{-1}(w). \een
Let
\ben \m=\g_{-1}\oplus\ldots\oplus \g_{-7}, \een 
where $gr_{-i}\cong\g_{-i}$. We have $\dim \m=\dim gr(w) =\dim T_w J^6=8$ and
\ben [\g_{-i},\g_{-j}]\subset\g_{-i-j}. \een

The additional piece of structure --- the distribution $\lspan\{\D\}$ --- is encoded in the following manner. One defines a weighted frame $z_w$ at $w\in J^6$ to be an isomorphism of graded Lie algebras $z_w\colon\m \to gr(w)$. The bundle of weighted frames $\mathcal{R} J^6$ is a principal bundle over $J^6$ with  the structure group $G_0(\m)$ being the group of all grading preserving algebra automorphisms of $\m$. 

The vector $\D_w$ at any $w$ belongs to $C_w$ and is complementary to the 1-dimensional subspace of $C_w$ which is vertical with respect to $J^6\to J^5$ and spanned by $\partial_{y_6}$. At the level of $\m$ it is reflected by a decomposition of $\g_{-1}$ into two 1-dimensional subspaces. These subspaces, call them $D$ and $V$ for short, are then encoded by reducing $\mathcal{R}J^6$ to a $G_0$-sub-bundle, where $G_0$ is the 2-dimensional subgroup of $G_0(\m)$ preserving the decomposition $\g_{-1}=D\oplus V$. 

However, the $G_0$-sub-bundle still is not the bundle $P$ where the connection $\Omega$ exists. In order to construct $P$ one needs to prolong the $G_0$-bundle. The procedure of prolongation is quite involved  and the reader is referred to the original paper \cite{Tan1}. The underling idea it is 
however simple. One aims to extend the $G_0$-bundle so that it is large enough to contain all the symmetries in the most symmetric homogeneous case.
From the example of the trivial ODE we know that the total space $P$ must be 11-dimensional, so one dimension is lacking. After the prolongation one obtains the desired $H\to P\to J^6$ with the structural group $H$ being a product of $G_0$ and the 1-dimensional prolongation, isomorphic to the group of triangular 2 by 2 matrices.

At the algebraic level the filtration is encoded by $\m$ and the full (prolonged) structural group $H$ is encoded by its algebra $\h=\g_0\oplus\g_1$. Since commutators $[\m,\h]$ are known from the construction, we obtain a graded algebra
\be\g=\m\oplus\h=\g_{-7}\oplus\g_{-6}\oplus\ldots\oplus\g_0\oplus\g_1=\gl\semi{.}\real^7. \ee
The Cartan connection $\Omega$ takes values in this algebra.

The next step is constructing the form $\Omega$ using the normality conditions of Tanaka and Morimoto. The normality conditions, which are certain linear constraints for the curvature, were originally introduced by E. Cartan in the context of conformal and projective geometries. The purpose was fixing ambiguity in the choice of Cartan connections and providing canonical connections for these geometries in a sense analogous to the Levi--Civita connection in Riemannian geometry. Later, these conditions were generalised to the case of the filtered manifolds. We discuss them below.

The connection 1-form at $p\in P$ is a vector space isomorphism $\Omega_p:T_pP\to\g$. 
We define $\V_p=\Omega_p^{-1}(\h)$ and $\H_p=\Omega_p^{-1}(\m)$, hence $T_pP= \V_p\oplus\H_p$. 
The curvature $K_p$ is then characterised by a tensor $\kappa_p\in\Hom(\wedge^2\m,\g)$ given by
 \be\label{e.kappa}\kappa_p(A,B)=K_p(\Omega^{-1}_p(A), \Omega^{-1}_p(B)),\quad
 A,B\in\m.\ee

In the space $\Hom(\wedge^2\m,\g)$ let us define
$\Hom^1(\wedge^2\m,\g)$ to be the space of all
$\alpha\in\Hom(\wedge^2\m,\g)$ fulfilling
\ben\alpha(\g_i,\g_j)\subset \g_{i+j+1}\oplus\ldots\oplus\g_k
\quad\text{for}\quad i,j<0.\een

The algebra $\g$ is equipped with the following complex \ben
\ldots\overset{\partial}{\longrightarrow}\Hom(\wedge^q\m,\g)
\overset{\partial}{\longrightarrow}\Hom(\wedge^{q+1}\m,\g)
\overset{\partial}{\longrightarrow}\ldots \een with
$\partial\colon\Hom(\wedge^q\m,\g)\to\Hom(\wedge^{q+1}\m,\g)$
given by \ben
\begin{aligned}
(\partial^*\alpha)&(A_1\w\ldots\w
A_{q+1})=\sum_i(-1)^{i+1}[A_i,\alpha(A_1\w\ldots\w
\hat{A}_i\w\ldots\w A_{q+1})]
 \\&+\sum_{i<j}(-1)^{i+j}\alpha([A_i,A_j]\w A_1\ldots\w
\hat{A}_i\w\ldots\w\hat{A}_j\w\ldots\w A_{q+1}),
\end{aligned}
\een where $\alpha\in\Hom(\wedge^q\m,\g)$ and
$A_1,\ldots\,A_{q+1}\in\m$. 

Consider a positive definite scalar product $(\cdot,\cdot)$ in $\g$
satisfying three conditions:
\begin{itemize}
 \item[i)] $(\g_i,\g_j)=0$ for $i\neq j$. \label{e.pair}
 \item[ii)] There exists a mapping $\tau\colon\h\to\g$ such that
\begin{align} 
\label{prop_of_in}
&\tau(\g_i)\subset\g_{-i}\qquad \text{for }i\geq 0, \\
&([A,X],Y)=(X,[\tau(A),Y])\qquad \text{for } X,Y\in\g,\,A\in\h. \nonumber
\end{align} 
\item[iii)] There exists a mapping $\tau_0\colon G_0\to G_0$ such that
\ben (aX,Y)=(X,\tau_0(a)Y)\qquad \text{for } X,Y\in\g,\,a\in G_0.  \een
\end{itemize}

This product extends to
$\Hom(\wedge^q\m,\g)$ through \ben
(\alpha,\beta)=\frac{1}{q!}\sum_{i_1,\ldots,i_q}
(\alpha(v_{i_1}\w\ldots\w v_{i_q}),\beta(v_{i_1}\w\ldots\w
v_{i_q})), \een where $\alpha,\beta\in\Hom(\wedge^q\m,\g)$ and
$(v_i)$ is any orthonormal basis of $\g$. Given $\partial$ and
$(\cdot,\cdot)$ the formal adjoint operator  \ben
\ldots\overset{\partial^*}{\longrightarrow}
\Hom(\wedge^{q+1}\m,\g)\overset{\partial^*}{\longrightarrow}
\Hom(\wedge^q\m,\g)\overset{\partial^*}{\longrightarrow}\ldots
\een is defined by \ben
(\partial^*\alpha,\beta)=(\alpha,\partial\beta).\een A normal connection is defined as follows.
\begin{defi}\label{def.norm}
A Cartan connection ${\Omega}$ is normal if 
$\kappa$ given by \eqref{e.kappa} satisfies
conditions \ben
\begin{aligned}
\text{i)}& &\qquad &\kappa\in\Hom^1(\wedge^2\m,\g),\\
\text{ii)}& &\qquad &\partial^*\kappa=0.
\end{aligned}
\een
\end{defi}
By a general result of Morimoto \cite{Mor2} (Theorem 2.3. and Proposition 2.10 in this reference), given an inner product satisfying the three properties (\ref{prop_of_in}) one can construct the normal Cartan connection which preserves the contact equivalence of the underlying ODEs. 
\section{Application to 7th order ODEs}
Define seven  1-forms on $J^6$ by:
\be
\label{e.om}
\omega^i=dy_{i-1}-y_idx, \quad \omega^7=dy_6-Fdx, \quad i=1, \dots, 6.
\ee
These forms encode the geometry of an ODE, and
in particular $C$ is annihilated by the ideal
$\lspan\{\omega^1,\omega^2,\ldots,\omega^6\}$ and the foliation by solutions is annihilated by
$\lspan\{\omega^1,\omega^2,\ldots,\omega^7\}$.
On $\tilde{P}$ there is the fundamental
$\real^8$-valued 1-form, whose components are denoted by
$\theta^1,\ldots,\theta^7$ and $\Gamma_+$. 
(The notation ${\Gamma}_+$ instead of
$\theta^8$ will be useful later on.)
One may introduce a coordinate system $(x,y,y_1,\ldots,y_6,u_1,u_2,\ldots,u_{36})$ compatible with the local trivialisation $\tilde{P}\cong
J^6\times\tilde{G}$ and such that locally 
\be \label{e.th}  \bma {\theta}^1 \\ {\theta}^2 \\
{\theta}^3 \\ {\theta}^4 \\ {\theta}^5 \\
{\theta}^6 \\ \theta^7  \\ {\Gamma}_+\ema = \bma u_1 & 0 & 0 & 0 & 0 & 0 & 0 & 0 \\
u_2 & u_3 & 0 & 0 & 0 & 0 & 0 & 0 \\
u_4 & u_5 & u_6 & 0 & 0 & 0 & 0 & 0 \\
u_7 & u_8 & u_9 & u_{10} & 0 & 0 & 0 & 0 \\
u_{11} & u_{12} & u_{13} & u_{14} & u_{15} & 0 & 0 & 0 \\
u_{16} & u_{17} & u_{18} & u_{19} & u_{20} & u_{21} & 0 & 0 \\
u_{22} & u_{23} & u_{24} & u_{25} & u_{26} & u_{27} & u_{28} & 0 \\
u_{29} & u_{30} & u_{31} & u_{32} & u_{33} & u_{34} & u_{35} &
u_{36} \ema \bma \omega^1 \\
\omega^2 \\ \omega^3 \\ \omega^4 \\
\omega^5 \\ \omega^6 \\ \omega^7 \\ \der x \ema. \ee 
The structural group $\tilde{G}$ is the group of the lower triangular 
matrices as above.

We choose a representation of $\gl\semi{.}\real^7$
and write down $\Omega$ in the following  matrix form
\be\label{e.con2}\Omega=\bma
-6\Gamma_0-6\Gamma_1 & 6\Gamma_+ & 0 & 0 & 0 & 0 & 0 & \theta^1\\\\
\Gamma_- & -4\Gamma_0-6\Gamma_1 & 5\Gamma_+ & 0 & 0 & 0 & 0 & \theta^2 \\\\
0 & 2\Gamma_- & -2\Gamma_0-6\Gamma_1 & 4\Gamma_+ & 0 & 0 & 0 & \theta^3 \\\\
0 & 0 & 3\Gamma_- & -6\Gamma_1 & 3\Gamma_+ & 0 & 0 & \theta^4\\\\
0 & 0 & 0 & 4\Gamma_- & 2\Gamma_0-6\Gamma_1 & 2\Gamma_+ & 0 & \theta^5\\\\
0 & 0 & 0 & 0 & 5\Gamma_- & 4\Gamma_0-6\Gamma_1 & \Gamma_+ & \theta^6\\\\
0 & 0 & 0 & 0 & 0 & 6\Gamma_- & 6\Gamma_0-6\Gamma_1 & \theta^7 \\\\
0 & 0 & 0 & 0 & 0 & 0 & 0 & 0 \ema. \ee
Here $\theta^1,\ldots\theta^7,\Gamma_+,\Gamma_0,\Gamma_1$ and $\Gamma_-$ are 1-forms on $P$.

Starting from this representation we construct a basis $(e_\mu)$, $\mu=1,\ldots, 11$ of $\gl\semi{.}\real^7$. 
To get the element $e_1$ we formally set $\theta^1=1$ and the remaining 1-forms equal to zero. 
All the remaining elements of the basis can be obtained in an analogous way, so that \eqref{e.con2} may be written as
\be\label{e.con1} \Omega=\sum^7_{i=1} \theta^ie_i+\Gamma_+e_8+\Gamma_0e_9+\Gamma_1e_{10}+\Gamma_-e_{11}. \ee
The basis satisfies
\begin{align*}
&\g_{-7}=\lspan\{e_1\}, & &\g_{-6}=\lspan\{e_2\}, & &\g_{-5}=\lspan\{e_3\},\\
&\g_{-4}=\lspan\{e_4\}, & &\g_{-3}=\lspan\{e_5\}, & &\g_{-2}=\lspan\{e_6\},\\
&\g_{-1}=\lspan\{e_7,e_8\},& &\g_{0}=\lspan\{e_9,e_{10}\}, & &\g_{1}=\lspan\{e_{11}\},
\end{align*}
and moreover
\begin{align*} &\gl=\lspan\{e_8,\ldots,e_{11}\}, && \real^7=\lspan\{e_1,\ldots,e_7\}. \end{align*}

To construct  $ P$ and ${\Omega}$ we need to
\begin{itemize}
\item[i)] Find a scalar product satisfying the conditions \eqref{prop_of_in}.
\item[ii)] Find formulae of $P\hookrightarrow \tilde{P}$ by expressing $u_4,\ldots,u_{36}$ as certain functions of
$u_1,u_2,u_3,x,y,y_1,\ldots,y_6$.  Then $(u_1,u_2,u_3,x,y,y_1,\ldots,y_6)$ is a local coordinate system in $P$ and the forms $\theta^1,\ldots,\theta^7,\Gamma_+$ of \eqref{e.con2} are given by the pull-back of \eqref{e.th}.
\item[iii)] Find formulae for $\Gamma_-$, $\Gamma_0$ and $\Gamma_1$.
\end{itemize}

We choose a scalar product on $\g$ so that the basis $(e_1,\ldots,e_{11})$ is orthogonal and
\begin{align*} &(e_1,e_1)=1, && (e_2,e_2)=6, &&(e_3,e_3)=15,\\
&(e_4,e_4)=20, &&(e_5,e_5)=15, &&(e_6,e_6)=6,\\
&(e_7,e_7)=1, &&(e_8,e_8)=1, &&(e_9,e_9)=2,\\
&(e_{10},e_{10})=1, &&(e_{11},e_{11})=1. && 
\end{align*}
The product satisfies the conditions \eqref{prop_of_in} if we set $\tau_0=id$, $\tau(e_9)=e_9$, $\tau(e_{10})=e_{10}$ and $\tau(e_{11})=e_8$. 

Both ii) and iii) are obtained from the horizontality condition \eqref{e.hor} and the normality conditions of Definition \ref{def.norm} with the scalar product as above. The 1-forms $\Gamma_0,\Gamma_1$ and $\Gamma_+$ on $P$ are a priori arbitrary
\ben
\Gamma_A=\,\sum_{j=1}^3 a^j_{~A}\,\der u_j+\sum_{i=1}^7b^i_{~A}\,\theta^i+b^+_{~A}\,\Gamma_+\qquad A=-,0,1.\\
\een
The functions $a$ and $b$ are arbitrary but sufficiently smooth  on $P$, so they depend on the jet coordinates and $u_1,u_2,u_3$.

The curvature \eqref{e.curv} becomes
\be\label{e.hor2} {K}=\sum_{\mu=1}^{11}\sum_{j=1}^7 K^\mu_{~8j}\Gamma_+\w\theta^j\otimes e_\mu+
\tfrac12\sum_{\mu=1}^{11}\sum_{i,j=1}^7 K^\mu_{~ij}\theta^i\w\theta^j\otimes e_\mu. \ee
and $K^\rho_{~\mu\nu}=-K^\rho_{~\nu\mu}$,
The terms proportional to $\Gamma_0,\Gamma_1$ and $\Gamma_-$ must be absent since ${K}$ is horizontal. 
This produces a set of first order differential equations for the functions $a$, which may be determined without ambiguity giving
\begin{align}\label{e.Gam_A}
\Gamma_0=&\,\frac12\frac{\der u_1}{u_1}-\frac12\frac{\der
u_3}{u_3}+\sum_i b^i_{~0}\theta^i+b^+_{~0}\Gamma_+,\\
\Gamma_1=&-\frac13\frac{\der u_1}{u_1}+\frac12\frac{\der
u_3}{u_3}+\sum_i b^i_{~1}\theta^i+b^+_{~1}\Gamma_+,\notag \\
\Gamma_-=&\,\frac{\der u_2}{u_1}+\frac{u_2\, \der u_3}{u_1
u_3}+\sum_i b^i_{~-}\theta^i+b^+_{~-}\Gamma_+. \notag
\end{align}

The tensor $\kappa$ is equal to
\ben
\kappa=\frac12\sum_{\mu=1}^{11}\sum_{i,j=1}^8 K^\mu_{~ij}e^i\w e^j\otimes e_\mu. \een
The condition
$\kappa\in\Hom^1(\wedge^2\m,\g)$ is equivalent to vanishing of the
following components of $K$.
\begin{align}\label{e.reg}
&K^{1,2,3,4,5}_{~67} &&K^{1,2,3,4}_{~57} && K^{1,2,3}_{~47}  \\
&K^{1,2}_{~37} &&K^{1}_{~27} && K^{1,2,3}_{~56} \notag \\
&K^{1,2}_{~46} &&K^{1}_{~36} && K^{1}_{~45} \notag \\
&K^{1,2,3,4,5,6}_{~87} &&K^{1,2,3,4,5}_{~86} && K^{1,2,3,4}_{~85} \notag \\
&K^{1,2,3}_{~84} &&K^{1,2}_{~83} && K^{1}_{~82}, \notag
\end{align}
where $K^{1,2}_{~37}$ is an abbreviation for $K^{1}_{~37}$,
$K^{2}_{~37}$ and so on. 

In order to evaluate the condition $\partial^*\kappa=0$ we introduce the notation $(e_\mu,e_\mu)=p_{\mu\mu}$ (no summation) and $[e_\mu,e_\nu]=c^\rho_{~\mu\nu}e_\rho$. The explicit form of $\partial^*\kappa=0$ is
\be\label{e.norm} 4\sum_{\nu=1}^{11}\sum_{j=1}^8 \frac{p_{\nu\nu}}{p_{ii}p_{jj}}K^\nu_{~ij}c^\nu_{~j\mu}
+\sum_{j,k=1}^8\frac{p_{\mu\mu}}{p_{jj}p_{kk}}K^\mu_{~jk}c^i_{~jk}=0, \ee
where $\mu=1,\ldots,11$ and $i=1,\ldots,8$.

We compute $K$ via \eqref{e.om}, \eqref{e.con2} and \eqref{e.Gam_A}. The
conditions \eqref{e.hor2} and \eqref{e.reg} become a set of easy
algebraic and differential equations for the functions
$u_4,\ldots,u_{36}$ and $b$. By solving these equations we
obtain $u_4,\ldots,u_{36}$ and $b$ as rational functions of
$u_1,u_2,u_3$ with coefficients given by  arbitrary functions of the
jet coordinates. After these substitutions the normality condition 
\eqref{e.norm} becomes a set of algebraic and differential
equations on the coefficients. The equations, although
complicated, are overdetermined and may be solved without
integration. It is enough to perform usual algebraic elimination of the
functions, provided it is done in an appropriate order. The
elimination also assures us that the solution --- the Cartan
connection --- is unique. We have therefore proved

\begin{prop}\label{th.equiv}
Given a 7th order ODE $y^{(7)}=F(x, y, y', \ldots,y^{(6)})$ one can construct
\begin{itemize}
\item[i)] A principal fibre bundle
$H\to P\to J^6$, where $H=\real\times(\real\ltimes\real)$.
\item[ii)] A Cartan connection ${\Omega}$ on
$P$ of type $(GL(2,\real)\ltimes\real^7,H)$.
\end{itemize}
Two ODEs
\ben y^{(7)}=F(x, y, y',\ldots,y^{(6)}) \een
and
\ben \bar{y}^{(7)}=\bar{F}(\bar{x}, \bar{y}, \bar{y}' ,\ldots, \bar{y}^{(6)}) \een
are locally contact
equivalent if and only if there exists a local bundle
diffeomorphism $\Phi\colon\bar{P}\to  P$
such that $\Phi^*\Omega=\,\bar{\Omega}$.
The connection is given by \eqref{e.con2}, where 

\begin{align*}
\theta^1&=u_1\omega^1, \notag \\
\theta^2&=u_2\omega^1+u_3\omega^2, \notag \\
\theta^3&= \frac{u_2^2}{u_1}\omega^1 + 2\frac{u_2u_3}{u_1}\omega^2 + \frac{u_3^2}{u_1}\Big(
\left (\frac{3}{14}(\D
F)_6-\frac{12}{35}F_5-\frac{13}{49}F^2_6\right
)\omega^1-\frac{2}{35}F_6\omega^2 +\frac65\omega^3\Big),\notag\\
\theta^4&=...&\\
\end{align*}
The explicit formulae for $\theta^i, i=4, 5, 7$ and $\Gamma_+, \Gamma_-, \Gamma_0, \Gamma_1$
are omitted since they are complicated and unilluminating. 
\end{prop}
\section{$\glg$ geometry from Cartan connection}\label{sec.gl2}
The manifold $P$ is endowed with two  structures of a principal bundle:
$H\to P\to J^6$ given by construction, and $\glg\to P\to M$
over the solution space which is generated by the connection $\Omega$. 
Let $X_\mu, \mu=1, \dots, 11$ denote the frame dual to the coframe
$(\theta^2,$ $\theta^2,$ $\theta^3,$ $\theta^4,$ $\theta^5,$ $\theta^6,$
$\theta^7,$ $\Gamma_+,$ $\Gamma_0,$ $\Gamma_1,$ $\Gamma_-)$ of \eqref{e.con2}. 
The curvature $K$ written in the form
\be\label{e.curv2}
\der\Omega=-\Omega\w\Omega+K, \ee
and split into scalar-valued equations reads

\begin{align}\label{e.dtheta}
\der\theta^1=&\,6(\Gamma_1+\Gamma_0)\w\theta^1-6\Gamma_+\w\theta^2+\tfrac12K^1_{~ij}\theta^i\w\theta^j,\\
\der\theta^2=&-\Gamma_-\w\theta^1+(6\Gamma_1+4\Gamma_0)\w\theta^2-5\Gamma_+\w\theta^3
+\tfrac12K^2_{~ij}\theta^i\w\theta^j,\nonumber\\
\der\theta^3=&-2\Gamma_-\w\theta^2+(6\Gamma_1+2\Gamma_0)\w\theta^3-4\Gamma_+\w\theta^4
+\tfrac12K^3_{~ij}\theta^i\w\theta^j\nonumber\\
&+K^3_{~18}\theta^1\w\Gamma_+,\nonumber\\
\der\theta^4=&-3\Gamma_-\w\theta^3+6\Gamma_1\w\theta^4-3\Gamma_+\w\theta^5
+\tfrac12K^4_{~ij}\theta^i\w\theta^j\nonumber\\
&+(K^4_{~18}\theta^1+K^4_{~28}\theta^2)\w\Gamma_+,\nonumber\\
\der\theta^5=&-4\Gamma_-\w\theta^4+(6\Gamma_1-2\Gamma_0)\w\theta^5-2\Gamma_+\w\theta^6
+\tfrac12K^5_{~ij}\theta^i\w\theta^j \nonumber \\
&+(K^5_{~18}\theta^1+K^5_{~28}\theta^2+K^5_{~38}\theta^3)\w\Gamma_+, \nonumber \\
\der\theta^6=&-5\Gamma_-\w\theta^5+(6\Gamma_1-4\Gamma_0)\w\theta^6-\Gamma_+\w\theta^7
+\tfrac12K^6_{~ij}\theta^i\w\theta^j\nonumber\\
&+(K^6_{~18}\theta^1+K^6_{~28}\theta^2+K^6_{~38}\theta^3+K^6_{~48}\theta^4)\w\Gamma_+,\nonumber\\
\der\theta^7=&-6\Gamma_-\w\theta^6+(6\Gamma_1-6\Gamma_0)\w\theta^7+\tfrac12K^7_{~ij}\theta^i\w\theta^j \nonumber\\
&+(K^7_{~18}\theta^1+K^7_{~28}\theta^2+K^7_{~38}\theta^3+K^7_{~48}\theta^4+K^7_{~58}\theta^5)\w\Gamma_+,\nonumber\\
\der\Gamma_+=&\,2\Gamma_0\w\Gamma_++\tfrac12K^8_{~ij}\theta^i\w\theta^j+K^8_{~i8}\theta^i\w\Gamma_+,\nonumber\\
\der\Gamma_0=&\,\Gamma_+\w\Gamma_-+\tfrac12K^9_{~ij}\theta^i\w\theta^j+K^9_{~i8}\theta^i\w\Gamma_+,\nonumber\\
\der\Gamma_1=&\tfrac12K^{10}_{~ij}\theta^i\w\theta^j+K^{10}_{~i8}\theta^i\w\Gamma_+,\nonumber\\
\der\Gamma_-=&-2\Gamma_0\w\Gamma_-+\tfrac12K^{11}_{~ij}\theta^i\w\theta^j+K^{11}_{~i8}\theta^i\w\Gamma_+.\nonumber
\end{align}

Since $J^6$ is a bundle over $M$ then so
is $P$ and Theorem \ref{th.equiv} together with Eqs. \eqref{e.om} guarantees 
that the fibres of the projection 
$P \to M$ are annihilated by the simple ideal $\lspan\{\theta^1,\ldots,\theta^7\}$.
The relation \eqref{e.dtheta} implies that this ideal is closed \ben
\der \theta^i\w\theta^1\w\ldots\w\theta^7=0 \qquad
\text{for}\qquad i=1,\ldots,7. \een 
It is  annihilated by an integrable distribution
$\lspan\{X_8,X_9,X_{10},X_{11}\}$ and the maximal integral
leaves of this distribution are locally the fibres of the projection $P \to M$. Moreover, by \eqref{e.dtheta} the
commutation relations of the vector fields are isomorphic to the
commutation relations of the algebra $\gl$. This allows us to {\em define} 
an action of $\glg$ on $P$ by defining $X_8,X_9,X_{10},X_{11}$ to be the associated fundamental vector fields.

\subsection{Existence of $\glg$ geometry.} Does the bundle $\glg\to
P\to M$ define a $\glg$ geometry on $M$? The answer to this
question is positive only if $P$ may be identified with a sub-bundle
of the frame bundle $ FM$. However, this may only be done
if the original ODE satisfies additional conditions. 
An object that turns $P$ into a sub-bundle of $F M$ is the canonical $\real^7$-valued
1-form. The $\real^7$-part of $\Omega$ --- the 1-forms
$\theta^1,\ldots,\theta^7$ arranged into a column --- is the natural candidate for it here, but one must
still check whether the canonical 1-form has the property  $R^*_g\theta=g^{-1}\theta$ under the actions of $\glg$ in $P$ and $\real^7$. 
In the language of the curvature $K$ this is equivalent to the horizontality with respect to the projection
$P\to M$. Since $K$ is already horizontal with respect to $P\to J^6$ we must only impose $K(X_8,\cdot)=0$ which amounts to
\ben K^\mu_{~\nu 8}=0,\qquad \mu,\nu=1,\ldots,11. \een
Due to algebraic and differential relations among the curvature components this condition may be further reduced to
\begin{align}\label{e.wunk} K^3_{~18}&=0, &  K^4_{~18}&=0, &
K^5_{~18}&=0, &  K^6_{~18}&=0, &  K^7_{~18}&=0 \end{align}
where
\[
K^{\alpha+2}_{~18}=\sum_{\beta=1}^\alpha {c^\alpha}_\beta W_\beta, \qquad \alpha=1, \dots, 5.
\]
The expressions $W_1, W_2, \dots, W_5$ are the W\"unschmann conditions discussed in the Introduction
and given by (Appendix A), and ${c^\alpha}_\beta$ are rational functions of $u_1, u_2, u_3$. The condition \eqref{e.wunk} is therefore equivalent to the vanishing of $W_\alpha$.
Simultaneous vanishing of these expressions is a property of a 7th order ODE invariant under contact transformations. It is also equivalent to the conditions for trivial linearizations obtained in \cite{Dou, Dou2}.

From now on we restrict our considerations to those ODEs which
satisfy all five conditions in Appendix A. Then the curvature contains no $\Gamma_+\w\theta^i$ terms, and the equations \eqref{e.dtheta}  may be written as the structural equations for a $\gl$-connection. We have proven
\begin{theo}\label{th.gl2}
Consider a 7th order ODE satisfying the conditions $W_\alpha=0$, $\alpha=1,\ldots, 5$. Then its solution space $M$ is equipped with a $\gl$ geometry. 
Let 
\be\label{e.gamma}\Gamma=\bma
-6\Gamma_0-6\Gamma_1 & 6\Gamma_+ & 0 & 0 & 0 & 0 & 0 \\\\
\Gamma_- & -4\Gamma_0-6\Gamma_1 & 5\Gamma_+ & 0 & 0 & 0 & 0 \\\\
0 & 2\Gamma_- & -2\Gamma_0-6\Gamma_1 & 4\Gamma_+ & 0 & 0 & 0 \\\\
0 & 0 & 3\Gamma_- & -6\Gamma_1 & 3\Gamma_+ & 0 & 0 \\\\
0 & 0 & 0 & 4\Gamma_- & 2\Gamma_0-6\Gamma_1 & 2\Gamma_+ & 0 \\\\
0 & 0 & 0 & 0 & 5\Gamma_- & 4\Gamma_0-6\Gamma_1 & \Gamma_+ \\\\
0 & 0 & 0 & 0 & 0 & 6\Gamma_- & 6\Gamma_0-6\Gamma_1 \ema. \ee 
be the $\gl$-part of the Cartan connection $\Omega$ of Proposition \ref{th.equiv} and $\theta=(\theta^i)$ be its $\real^7$ part. 
Then $\Gamma$ is a $\gl$ linear connection on $P$ compatible with the $\glg$ geometry and the equations \eqref{e.curv}, \eqref{e.dtheta} read
\begin{align} 
&\der\theta^i+\Gamma^i_{~j}\wedge\theta^j= \frac12
T^i_{~kl}\theta^k\w\theta^l,\qquad i,j=1,\ldots, 7, \label{e.str1} \\ 
&\der\Gamma^{i}_{~j}+\Gamma^i_{~k}\w\Gamma^k_{~j}=\frac12
R^i_{~jlm}\theta^l\w\theta^m, \label{e.str2}
\end{align} 
where $T$ and $R$ are the torsion and curvature of $\Gamma$ respectively.
\end{theo}

We will construct two tensor fields on $M$ preserved by the $\glg$ geometry: the conformal metric $g$ and the conformal 3-form $\phi$. This is done as follows. The action of $\gl$ on $\real^7$ is given by the matrix representation \eqref{e.gamma}, in particular it defines two the conformal classes of tensors represented by
 $g\in S^2\real^{7*}$ and $\phi\in\Lambda^3\real^{7*}$. Next we transport these tensors to $T^*P$. The connection $\Omega$ gives the identification $e^i\leftrightarrow\theta^i$, $i=1,\ldots,7$ where $(e^i)$ is dual of the basis $(e_i)$  of \eqref{e.con2} and \eqref{e.con1}. By this identification we get the tensor fields on $P$:
\ben g=\theta^{{1}}\theta^{{7}}-6\,\theta^{{2}}\theta^{{6}}+15\,\theta^{{3}}
\theta^{{5}}-10\,(\theta^4)^2\een
and
\begin{align} \phi=&3\,\theta^{{2}}\w\theta^{{3}}\w\theta^{{7}}
-6\,\theta^{{2}}\w\theta^{{4}}\w\theta^{{6}} -\theta^{{1}}\w
\theta^{{4}}\w\theta^{{7}}
+3\,\theta^{{1}}\w\theta^{{5}}\w\theta^{{6}}+
15\,\theta^{{3}}\w\theta^{{4}}\w\theta^{{5}}. \label{e.fi}
\end{align}
Finally, we project these fields to conformal fields on $M$. This projection is well-defined because $g$ and $\phi$ satisfy two conditions: i) the vertical directions of $P \to M$ are degenerate for $g$ and $\phi$, and ii) the vertical directions are conformal symmetries of $g$ and $\phi$, that is 
\begin{align*} {\mathcal L}_{X_8}g=&0 & {\mathcal L}_{X_9}{g}&=0,& {\mathcal
L}_{X_{10}}{g}&=12\,{g},& {\mathcal L}_{X_{11}}{g}&=0, \\
{\mathcal L}_{X_8}\phi=&0 & {\mathcal L}_{X_9}{\phi}&=0,& {\mathcal
L}_{X_{10}}{\phi}&=18\,{\phi},& {\mathcal L}_{X_{11}}{\phi}&=0.
\end{align*}
It is worth noting that ${\mathcal L}_{X_8}g={\mathcal L}_{X_8}\phi=0$ are equivalent to conditions listed in Appendix A.
The conformal fields on $M$ will be also denoted by $g$ and $\phi$ --- on solutions to the
ODE they coincide with  (\ref{7Dconf}) and (\ref{three_form}) respectively. 

The following fact is an immediate consequence of Theorem \ref{th.gl2}.
\begin{prop}\label{th.cov}
Let $\nabla$ denote the covariant derivative on $M$ associated to $\Gamma$. We have
\begin{align*}
&\nabla_Xg=-A(X)g,\\
&\nabla_X\phi=-\tfrac32A(X)\phi,
\end{align*}
where the 1-form $A$ is proportional to the trace of the
connection matrix: \ben A=\tfrac27\sum_j \Gamma^j_{~j} = \sum_{i,j}\langle\nabla_i X_j, \xi^j\rangle \xi^i,\een for
any frame $(X_i)$ and the dual coframe $(\xi^i)$ such that $g=g_{ij}\xi^i\otimes \xi^j$ with constant $g_{ij}$.
\end{prop}

Of course, $g$ and $\phi$ do not reduce $GL(7,\real)$ to $\glg$, since their conformal stabilisers are $CO(3,4)$ and $\real^+\times G^{split}_2$ respectively. The object whose conformal stabiliser is precisely the irreducible $\glg$ is a certain totally symmetric 4-tensor $\Upsilon_{ijkl}$, which is however irrelevant in our approach.
\subsection{Torsion}
In this section we consider only those ODEs which admit  the $\glg$ geometry on the solution space. First  we shall characterise the torsion $T$ of $\Gamma$.
Let $V^k$ denote the $k$-dimensional irreducible representation of $\glg$ as before. 
Torsion of any $\gl$-connection at $p\in
P$ belongs to the representation $\Lambda^2 V^{7*}\otimes V^7$ which
decomposes as  \ben \Lambda^2 V^{7*}\otimes V^7=V^1\oplus
V^3 \oplus 3 V^5 \oplus 3V^7\oplus 3V^9\oplus 2V^{11}\oplus
2V^{13}\oplus V^{15}\oplus V^{17}. \een 
\begin{prop}
\label{prop_tor}
The only non-vanishing components of 
the torsion $T$ of the connection $\Gamma$ in Theorem \ref{th.gl2} are 
in the 1-dimensional, the
3-dimensional, and a fixed 5-dimensional representation in the
above decomposition. \[ T=T_1+T_3+T_5.\]
Explicit form of $T$ in \eqref{e.str1} is given in Appendix B, where $\lambda$ spans $T_1$; $a_1,a_2,a_3$ span $T_3$,
and $b_1,b_2,b_3,b_4,b_5$ span $T_5$.
\end{prop}
{\bf Proof.}
To prove the formula of Appendix B we use  Proposition \ref{th.equiv} to

explicitly calculate \eqref{e.str1}. Next we check that $T$ only occupies the irreducible representations as above.
\koniec

We are now ready to prove Theorems  \ref{th.g2} and \ref{th.examples}.

{\bf Proof of Theorem \ref{th.g2}.}
Using \eqref{e.fi}, \eqref{e.str1} and \eqref{FG_deco} we calculate $\lambda$, $\tau_2$, $\tau_3$ and $\Theta$ in terms of
the torsion coefficients and the forms $\theta^i$. We find that 
$\tau_2=0$ if $T_3=0$, $\tau_3=0$ if $T_5=0$, and also $\Theta=24\Gamma_1$. Next we
calculate explicitly $\lambda$, $a_i$ and $b_i$ using formulae for $\Omega$ given in Theorem \ref{th.equiv}. Since the components $T_3$ and $T_5$ lie in irreducible representations they vanish iff any
of the components $a_i$ or $b_i$ vanishes. In the Theorem we gave
the simplest ones.\koniec

{\bf Proof of Theorem  \ref{th.examples}.}
In order to prove this result we need to extensively use the
Bianchi identities. First, we suppose that $\tau_3= 0$ which
is equivalent to vanishing of $b_1, b_2, \dots, b_5$. Then from $\der ^2
\theta^i= 0$ we find that either i) $a_i= 0$
(equivalently $\tau_2= 0$) or ii) $\lambda= 0$.

Suppose i). Then the torsion is reduced to $\lambda$ and it makes
all the curvature except the Ricci scalar vanish. In particular the Lee form $\Theta=24\Gamma_1$ is closed. Therefore there exists  a conformal gauge in which locally $\Theta= 0$ and
$\lambda=const$, and which defines a 10-dimensional sub-bundle $P'$
of $P$. Eqs. \eqref{e.str1}, \eqref{e.str2} pulled-back to $P'$
become the structural equations of $SO(3,2)$ while the integrable
distribution on $P'$ annihilated by $\theta^1,\ldots,\theta^7$
defines the action of $SO(2,1)$ on $P'$, which is vertical w.r.t.
$P'\to M$. This also means that the maximal symmetry group of an
underlying ODE is $SO(3,2)$. We find the ODE of point 2. by
integration of the conditions from Appendix A and the conditions of
Theorem \ref{th.g2}.

Suppose ii). Lengthy but straightforward
calculations show that the condition $\lambda=0$ specifies curvature in Eqs. 
\eqref{e.str2}; all torsion and curvature coefficients and their
coframe derivatives are polynomials of $a_1,a_2$ and $a_3$, which
span $T_3$. Again, we have $\der \Theta= 0$. Since $T_3$
belongs to the 3-dimensional representation $V^3$ of $\glg$ we may
classify it by the orbit is sweeps out. $T_3$ is a tensor field on
$P$, which is a $\gl$ bundle over $M$. If we fix $x\in M$ and
sweep out the fibre $P_x$ then $T_3$ at points $p\in P_x$ sweeps a
$\glg$-orbit in $V^3$. These orbits are labelled by the sign of
$\langle\cdot,\cdot\rangle$, the conformal product in $V^3$
preserved by $\glg$. The case $\langle T_3,T_3 \rangle
= 0$ is forbidden by the Bianchi identities. The only
remaining possibilities are $\langle T_3,T_3\rangle >0$ and
$\langle T_3,T_3\rangle<0$. The ODE (\ref{cusp_ODE}) generates both
cases in two disjoint areas of $M$, depending on the sign of
$\langle T_3, T_3 \rangle= const\cdot(5y^{(6)} y^{(4)}-6(y^{(5)})^2)$.
\koniec
\section*{Appendix A}
\begin{appendix}
The five Wunschmann conditions for the 7th order ODE
\begin{eqnarray*}
&&W_1=245 \D ^2F_6 - 245 \D F_5 + 98 F_4 - 210 \D F_6 F_6 + 70 F_5 F_6 + 20F_6^3\\&&\\
&&W_2=6860 \D ^2F_5 - 10976 \D F_4 + 6615 (\D F_6)^2 + 6860 F_3 - 8330 \D F_6 F_5 +\\&&
    1715 F_5^2 - 1960 \D F_5 F_6 + 1568 F_4 F_6 - 1890 \D F_6 F_6^2 +
    1190 F_5 F_6^2 + 135 F_6^4\\&&\\
&&W_3=9604 \D ^2F_4 - 24010 \D F_3 + 15435 \D F_5 \D F_6 + 24010 F_2 - 14749 \D F_6
  F_4 -
\\&&
    5145 \D F_5 F_5 + 4459 F_4 F_5 - 2744 \D F_4 F_6 + 6615 (\D F_6)^2 F_6 +
    3430 F_3 F_6 -
\\&& 6615 \D F_6 F_5 F_6 + 1470 F_5^2 F_6 - 2205 \D F_5 F_6^2 +
    2107 F_4 F_6^2 -\\&& 1890 \D F_6 F_6^3 + 945 F_5 F_6^3 + 135
    F_6^5\\&&\\
&&W_4=336140 \D ^2F_3 - 1344560 \D F_2 + 180075 (\D F_5)^2 + 432180 \D F_4 \D F_6 +\\&&
    2352980 F_1 - 624260 \D F_6 F_3 - 216090 \D F_5 F_4 + 64827 F_4^2 - \\&&
    144060 \D F_4 F_5 + 154350 (\D F_6)^2 F_5 + 192080 F_3 F_5 -
    102900 \D F_6 F_5^2 +\\&& 17150 F_5^3 - 96040 \D F_3F_6 +
    308700 \D F_5 \D F_6 F_6 + 192080 F_2 F_6 -\\&& 246960 \D F_6 F_4 F_6 -
    154350 \D F_5 F_5 F_6 + 113190 F_4 F_5 F_6 - 61740 \D F_4 F_6^2 +\\&&
    132300 (\D F_6)^2 F_6^2 + 89180 F_3 F_6^2 - 176400 \D F_6 F_5 F_6^2 +
    47775 F_5^2 F_6^2 -\\&& 44100 \D F_5 F_6^3 + 35280 F_4F_6^3 -
    37800 \D F_6F_6^4 + 22050 F_5 F_6^4 + 2700 F_6^6\\&&\\
&&W_5=2352980 \D ^2F_2 - 16470860 \D F_1 + 1512630 \D F_4 \D F_5 +
      2268945 \D F_3 \D F_6 - \\&&5126135 \D F_6 F_2 - 1512630 \D F_5 F_3 -
      907578 \D F_4 F_4 + 648270 (\D F_6)^2 F_4 +\\&& 907578 F_3 F_4 -
      756315 \D F_3 F_5 + 1080450 \D F_5 \D F_6 F_5 + 1596665 F_2 F_5 - \\&&
      1080450 \D F_6 F_4F_5 - 360150 \D F_5 F_5^2 + 288120 F_4 F_5^2 -
      672280 \D F_2 F_6 +\\&& 540225 (\D F_5)^2F_6 + 1296540 \D F_4 \D F_6 F_6 +
      2352980 F_1 F_6 -\\&& 1620675 \D F_6 F_3 F_6 - 864360 \D F_5 F_4 F_6 +
      324135 F_4^2 F_6 - 648270 \D F_4 F_5 F_6 +\\&& 926100 (\D F_6)^2 F_5 F_6 +
      756315 F_3 F_5 F_6 - 771750 \D F_6 F_5^2 F_6 + 154350 F_5^3 F_6 -\\&&
      324135 \D F_3 F_6^2 + 926100 \D F_5 \D F_6 F_6^2 + 732305 F_2 F_6^2 -
      926100 \D F_6 F_4 F_6^2 -\\&& 617400 \D F_5 F_5 F_6^2 +
      524790 F_4 F_5 F_6^2 - 185220 \D F_4 F_6^3 + 396900 (\D F_6)^2 F_6^3
  +
\\&&
      231525 F_3 F_6^3 - 661500 \D F_6 F_5 F_6^3 + 209475 F_5^2 F_6^3 -
      132300 \D F_5 F_6^4 +\\&& 119070 F_4 F_6^4 - 113400 \D F_6 F_6^5 +
      75600 F_5 F_6^5 + 8100 F_6^7 + 65883440 F_0.
\end{eqnarray*}
\end{appendix}
\section*{Appendix B}
\begin{appendix}
The torsion components in Proposition \ref{prop_tor}
\begin{align*}
T^{{1}}=&{\tfrac {55}{18}}{ b_1} \theta^{{1}}\w\theta^{{2}}
+{\tfrac {55}{9}}{ b_4} \theta^{{1}}\w\theta^{{3}} +
\left({\tfrac{55}{18}}{ b_3}-\tfrac{10}{3}\lambda-3{ a_3} \right)
\theta^{{1}}\w\theta^{{4}}\\
&+ \left( -{ \tfrac {55}{9}}{ b_5}+\tfrac{3}{2}{ a_2} \right)
\theta^{{1}}\w\theta^{{5}}-{\tfrac {77}{36}}{ b_2}
\theta^{{1}}\w\theta^{{6}} +\left( -{\tfrac {55}{2}}{
b_3}+10\lambda+9{ a_3} \right) \theta^{{2}}\w\theta^{{3}}\\
&+ \left( {\tfrac {55}{3}}{ b_5}-3{ a_2} \right)
\theta^{{2}}\w\theta^{{4}} +{\tfrac{55}{12}} { b_2}
\theta^{{2}}\w\theta^{{5}} ,\\
T^{{2}}=&{\tfrac {55}{36}}{ b_1}
\theta^{{1}}\w\theta^{{3}} +
\left( {\tfrac {275}{54}}{ b_4}+\tfrac12{ a_1} \right)
\theta^{{ 1}}\w\theta^{{4}}
+ \left( -{\tfrac {55}{36}}{ b_3}-\tfrac53 \lambda-{ a_3}\right)
\theta^{{1}}\w \theta^{{5}} \\
&+\left(-{\tfrac {11}{18}}{ b_5}+\tfrac12{ a_2 } \right)
\theta^{{1}}\w\theta^{{6}}
-{\tfrac {77}{216}}{ b_2}\theta^{{1}}\w\theta^{{7}}
+ \left( -{\tfrac{55}{18}}{ b_4} -\tfrac32{ a_1} \right)
\theta^{{2}}\w\theta^{{3}} \\
&+\left( -{\tfrac {55}{9}}{ b_3}+2{ a_3}+\tfrac{10}{3}\lambda
\right) \theta^{{2}}\w \theta^{{4}}-{\tfrac {11}{18}}{ b_2}
\theta^{{2}}\w\theta^{{6}} +\left( {\tfrac {275}{18}} {
b_5}-\tfrac52{ a_2} \right) \theta^ {{3}}\w\theta^{{4}}\\
&+{\tfrac {275}{72}}{ b_2}
\theta^{{3}}\w\theta^{{5}} ,\\
T^{{3}}=&{\tfrac{11}{54}}{ b_1}
\theta^{{1}}\w\theta^{{4}} +
\left( {\tfrac {22}{9}}{ b_4}+\tfrac12{ a_1}\right)
\theta^{{1}}\w\theta^{{5}} +
\left( -{\tfrac {44}{45}}{ b_3}-\tfrac15{ a_3}-\tfrac23\lambda \right)
\theta^{{1}}\w\theta^{{6}} \\
&+\left({\tfrac {22}{135}}{ b_5}+\tfrac{1}{10}{ a_2} \right)
\theta^{{1}}\w\theta^{{7}} +
{ \tfrac {22}{9}}{ b_1}
\theta^{{2}}\w \theta^{{3}} +
\left( {\tfrac{11}{9}}{ b_4}-{ a_1}\right)
\theta^{{2}}\w\theta^{{4}}\\
&+\left( -{\tfrac {11}{45}}{ b_5}+\tfrac35{ a_2 } \right)
\theta^{{2}}\w\theta^{{6}}
-{\tfrac {11}{20}}{ b_2}
\theta^{{2}}\w\theta^{{7}} +
\left( {\tfrac {55}{18}}{ b_3}+ \tfrac{10}{3}\lambda+{ a_3} \right)
\theta^{{3}}\w\theta^{{4}} \\
&+\left( {\tfrac {55}{9}}{ b_5}-\tfrac32{ a_2} \right) 
\theta^{{3}}\w\theta^{{5}}
+{\tfrac {11}{12}}{ b_2} \theta^{{3}}\w\theta^{{6}}
+{\tfrac {55}{18}}{ b_2} \theta^{{4}}\w\theta^{{5}}
-\tfrac{11}{2}{ b_3} \theta^{{2}}\w\theta^{{5}},\\
T^{{4}}=&-{ \tfrac {11}{24}}{ b_1} \theta^{{1}}\w\theta^{{5}} +
\left( {\tfrac {11}{15}}{ b_4}+\tfrac{3}{10}{ a_1 } \right)
\theta^{{1}}\w\theta^{{6}} + \left( -{\tfrac {11}{90}}{
b_3}-\tfrac16\lambda \right) \theta^{{1}}\w\theta^{{7}}\\
&+{\tfrac {11}{6}}{ b_4} \theta^{{2}}\w\theta^{{5}} + \left(
-{\tfrac {22}{5}} { b_3}-\lambda \right) \theta^{{2}}\w
\theta^{{6}} + \left( {\tfrac {11}{15}}{ b_5}+\tfrac{3}{10}{ a_2 }
\right) \theta^{{2}}\w\theta^{{7}}\\
&+ \left( {\tfrac {55}{18}}{b_4}-\tfrac32{ a_1} \right)
\theta^{{3}}\w\theta^{{4}}+\tfrac52 \lambda
\theta^{{3}}\w\theta^{{5}} +{\tfrac {11}{6}}{ b_5}
\theta^{{3}}\w\theta^{{6}} -{\tfrac {11}{24}}{ b_2}
\theta^{{3}}\w\theta^{{7}}\\
&+ \left( { \tfrac{55}{18}}{ b_5}-\tfrac32{ a_2} \right)
\theta^{{4}}\w\theta^{{5}}
+{\tfrac {22}{9}}{ b_2}\theta^{{4}}\w\theta^{{6}}
+{\tfrac {22}{9}}{ b_1} \theta^{{2}}\w \theta^{{4}}, \\
T^{{5 }}=&-{\tfrac {11}{20}}{ b_1} \theta^{{1}}\w\theta^{{6}} +
\left( {\tfrac {22}{135}}{ b_4}+\tfrac{1}{10}{  a_1} \right)
\theta^{{1}}\w\theta^{{7}} +{\tfrac {11}{12}}{ b_1}
\theta^{{2}}\w\theta^{{5}} \\
&+\left( -{\tfrac {11}{45}}{ b_4} +\tfrac35{ a_1} \right)
\theta^{{2}}\w \theta^{{6}}+ \left( -{\tfrac {44}{45}}{
b_3}-\tfrac23\lambda +\tfrac15{ a_3} \right) \theta^{{2}}\w
\theta^{{7}} \\
&+{\tfrac{55}{18}}{ b_1} \theta^{{3}}\w\theta^{{4}} +\left(
{\tfrac {55}{9}} { b_4}-\tfrac32{ a_1} \right)
\theta^{{ 3}}\w\theta^{{5}}-\tfrac{11}{2}{ b_3} \theta^{{3}}\w\theta^{{6}}\\
&+\left( {\tfrac {22}{9}}{ b_5}+\tfrac12{ a_2} \right) \theta^{{3
}}\w\theta^{{7}} + \left( {\tfrac {55}{18}}{ b_3}+\tfrac{10}{3}
\lambda-{ a_3}\right) \theta^{{4}}\w \theta^{{5}}\\
&+ \left( {\tfrac{11}{9}}{ b_5}-{ a_2} \right)
\theta^{{4}}\w\theta^{{6}}+{\tfrac {11}{54}}{ b_2}
\theta^{{4}}\w\theta^{{7}}
+{\tfrac {22}{9}}{ b_2} \theta^{{5}}\w\theta^{{6}} ,\\
T^{{6}}=&-{ \tfrac {77}{216}}{ b_1}
\theta^{{1}}\w\theta^{{7}}
-{\tfrac {11}{18}}{ b_1} \theta^{{2}}\w\theta^{{6}}
+ \left( -{\tfrac { 11}{18}}{ b_4}+\tfrac12{ a_1} \right)
\theta^{{2}}\w\theta^{{7}}\\
&+\left( -{\tfrac {55}{36}}{ b_3}-\tfrac53\lambda+{ a_3} \right)
\theta^{{3}}\w\theta^{{7}} + \left( {\tfrac {275}{18}}{
b_4}-\tfrac52{ a_1} \right) \theta^{{4}}\w\theta^{{5}}\\
&+ \left(-{\tfrac { 55}{9}}{ b_3}+\tfrac{10}{3}\lambda-2{ a_3}
\right) \theta^{{4}}\w\theta^{{6}}+\left( {\tfrac { 275}{54}}{
b_5}+\tfrac12{ a_2} \right) \theta^{{4}}\w\theta^{{7}}\\
&+ \left( -{\tfrac {55}{18}}{ b_5}-\tfrac32{ a_2} \right)
\theta^{{5}}\w\theta^{{6}}
+{\tfrac {55}{36}}{ b_2}\theta^{{5}}\w\theta^{{7}}
+{\tfrac {275}{72}}{ b_1} \theta^{{3}}\w\theta^{{5}} ,\\
T^{{7}}=&-{\tfrac { 77}{36}}{ b_1} \theta^{{2}}\w\theta^{{7}}
+{\tfrac {55}{12}}{ b_1}\theta^{{3}}\w\theta^{{6}} + \left(
-{\tfrac {55}{9}}{ b_4}+\tfrac32{ a_1} \right) \theta^{{3
}}\w\theta^{{7}}\\
&+ \left( {\tfrac {55}{3}}{ b_4}-3{ a_1}\right)
\theta^{{4}}\w\theta^{{6}} + \left( {\tfrac {55}{18}}{
b_3}-\tfrac{10}{3}\lambda+3{ a_3}\right)
\theta^{{4}}\w\theta^{{7}}\\
&+\left( -{\tfrac {55}{2}}{
b_3}+10\lambda-9{ a_3}\right) \theta^{{5}}\w\theta^{{6}} +{\tfrac
{55}{9}}{ b_5} \theta^{{5}}\w\theta^{{7}} +{\tfrac {55}{18}}{
b_2}\theta^{{6}}\w\theta^{{7}}.
\end{align*}
\end{appendix}


\begin{thebibliography}{jafsdl}
\frenchspacing 


\bibitem{Nur_Bob} Bobie\'nski, M. and Nurowski, P. (2007)
Irreducible SO(3) geometries in dimension five,  J. reine angew Math. 
{\bf 605}. 

\bibitem{B87} Bryant, R. L. (1987) 
Metrics with exceptional holonomy.  Ann. of Math. (2)  126,  
no. 3, 525--576.  
\bibitem{B91} Bryant, R. L. (1991) Two exotic holonomies in dimension four,
path geometries, and twistor theory.
Proc. Symp. Pure. Maths. Vol. 53, 33--88. 
\bibitem{B03} Bryant, R. L. (2003)  Some remarks on $G_2$-structures,
{\tt arXiv:math/0305124v4}.
\bibitem{ivanov} Cleyton, R. \& Ivanov, S. (2007)
On the geometry of closed $G_2$-structures. Comm. Math. Phys. {\bf 270}, 53--67.
\bibitem{Dou} Doubrov, B. (2001)
Contact trivialization of ordinary differential equations. 
Differential geometry and its applications, 73--84.
\bibitem{Dou2} Doubrov, B. (2008)
Generalized Wilczynski invariants for non-linear ordinary differential equations,
The IMA Volumes in Mathematics and its Applications {\bf 144}, 25-40.
\bibitem{DKM}  Doubrov, B., Komrakov, B. \& Morimoto, T. (1999)
Equivalence of holonomic differential equations, {\it Lobachevskij
Journal of Mathematics} {\bf 3} 39.

\bibitem{DD} Doubrov, B. \& Dunajski, M. (2011) Co--calibrated $G_2$ structure from cuspidal cubics. {\tt arXiv:1107.2813}. To appear
in Annals of Global Analysis and Geometry.

\bibitem{DT06} Dunajski, M. \& Tod, K. P. (2006) 
{Paraconformal geometry of  $n$th order ODEs, and exotic
holonomy in dimension four.}  J. Geom. Phys. 56, 1790-1809.
\bibitem{DS10} Dunajski, M. \& Sokolov, V. V. (2010) On 7th order ODE with submaximal symmetry. Preprint {\tt arXiv:1002.1620}.
\bibitem{FG_paper} Fernandez, M.,  Gray, A. (1982)
Riemannian manifolds with structure group $G_2$, Ann. Mat.
Pura Appl. {\bf 132}  19-45.
\bibitem{Nur} Godli\'nski, M. \& Nurowski, P. 
(2010) GL$(2, \R)$ geometry of ODEs, J. Geom. Phys. {\bf 60}, 991-1027.
\bibitem{Grace_Young} 
Grace, J. H. and Young, Alfred (1903). {\it The algebra of invariants.} 
Cambridge: Cambridge University Press.
\bibitem{hitchin} Hitchin, N.\ (1982) Complex manifolds and Einstein's
equations, in {\em Twistor Geometry and Non-Linear systems}, Springer 
LNM 970, Doebner, H.\ \& Palev. T.
\bibitem{Kob} Kobayashi, S. {\it Transformation Groups in Differential Geometry}.
\bibitem{kodaira} Kodaira, K. (1963)
On stability of compact submanifolds of complex manifolds,
 Am. J. Math.  {\bf 85}, 79-94.
\bibitem{Mor1} Morimoto, T. (1993) Geometric structures on filtred manifolds,
{\it Hokkaido Math. J.} {\bf 22}  263.
\bibitem{Mor2} Morimoto, T. Lie Algebras, Geometric Structures
and Differential Equations on Filtred Manifolds, in: {\it Advanced
Studies in Pure Mathematics} 37, 2002, Lie Groups, Geometric
Structures and Differential Equations --- One Hundred Years after
Sophus Lie, 205-252.

\bibitem{Nurowski} Nurowski, P. (2005)
Differential equations and conformal structures.  J. Geom. Phys.  {\bf 55}, 
19--49.

\bibitem{Nur_four} Nurowski, P. (2009) 
Comment on GL(2,R) geometry of 4th order ODEs,
J. Geom. Phys. {\bf 59}. 

\bibitem{olver} Olver, P. J. (1995) {\em Equivalence, invariants, 
and symmetry.} CUP, Cambridge.

\bibitem{PR} Penrose, R. \& Rindler, W.
 (1987, 1988) 
{\em Spinors and space-time.  Two-spinor calculus and relativistic fields. }
Cambridge Monographs on Mathematical Physics. Cambridge University Press, 
Cambridge.
\bibitem{sokolov}  Sokolov, V. V. (1988) 
Symmetries of evolution equations. Russian Math. Surveys  {\bf 43} 
165--204.
\bibitem{Tan1} Tanaka, N. (1970) On differential systems, graded Lie algebras and pseudo-groups, {\it J. Math. Kyoto Univ} {\bf 10}  1.
\bibitem{Tan2} Tanaka, N. (1979) On the equivalence probles associated
with simple graded Lie algebras, {\it Hokkaido Math. J.} {\bf 8}
23.
\bibitem{morozov}  Vinberg, E. B., Gorbatsevich, V. V. \& Onishchik, A. L. (1994)
{\em Structure of Lie groups and Lie algebras. 
Lie groups and Lie algebras, III.} 
Encyclopaedia of Mathematical Sciences, {\bf 41}. Springer-Verlag, Berlin.

\bibitem{alg_geom_book} Walker, R (1950){\em Algebraic Curves}, Princeton.

\bibitem{Wilczynski} Wilczynski, E. J. (1905) {\em Projective differential geometry of curves and ruled surfaces}, Leipzig, Teubner.
\bibitem{wun} W\"unschmann, K.W. (1905) {\em Inaug. Dissert.},  
Teubner, Leipzig.

\end{thebibliography}
\end{document}